\documentclass[oneside,11pt]{article}
\usepackage[top=1.25in, bottom=1.25in, left=0.8in, right=0.8in]{geometry}
\usepackage{amssymb}
\usepackage{amsthm}
\usepackage{graphicx}
\usepackage{setspace}
\usepackage{titling}
\usepackage[utf8]{inputenc}
\usepackage[english]{babel}
\usepackage{mathrsfs}
\usepackage{amsfonts}
\usepackage{amsmath}
\usepackage{float}
\usepackage{color}
\usepackage{dsfont}
\usepackage[toc,page]{appendix}
\usepackage{eqnarray}
\usepackage{bbm}
\usepackage{mathtools}
\usepackage{verbatim}

\usepackage{accents}

\newtheorem{theorem}{Theorem}

\def\bW{{\boldsymbol W}}
\def\bX{{\boldsymbol X}}
\def\bY{{\boldsymbol Y}}

\def\BB{\mathbb B}
\def\CC{\mathbb C}
\def\DD{\mathbb D}
\def\EE{\mathbb E}
\def\FF{\mathbb F}
\def\LL{\mathbb L}
\def\PP{\mathbb P}
\def\RR{\mathbb R}
\def\SS{\mathbb S}
\def\XX{\mathbb X}
\def\YY{\mathbb Y}
\def\ZZ{\mathbb Z}
\def\balpha{\boldsymbol{\alpha}}

\def\bmu{\boldsymbol{\mu}}
\def\cC{\mathcal C}
\def\cF{\mathcal F}
\def\cL{\mathcal L}
\def\cP{\mathcal P}

\begin{document}

\title{An Alternative Approach to Mean Field Game with Major and Minor Players, and Applications to Herders Impacts}
\author{Rene Carmona \and Peiqi Wang}
\maketitle

\begin{abstract}
The goal of the paper is to introduce a formulation of the mean field game with major and minor players as a fixed point on a space of controls.
This approach emphasizes naturally the role played by McKean-Vlasov dynamics in some of the players optimization problems. 
We apply this approach to linear quadratic models for which we recover the existing solutions for open loop equilibria, and we show that we
can also provide solutions for closed loop versions of the game. Finally, we implement numerically our theoretical results on a simple model of flocking.
\end{abstract}

\section{Introduction}
Mean field games with major and minor players were introduced with the specific intent to extend the realm of applications of the original mean field game paradigm to realistic models for which subgroups of players do not grow in size and as a result, their influence on the remaining population of players, does not disappear in the asymptotic regime of large games. While this generalization captures new potential applications, it \emph{raises the technological bar} in terms of the sophistication of the tools to be used in order to come up with solutions, bringing these models up to par with mean field games with common noise. See for example the monograph \cite{BensoussanChauYam} or the last chapter of \cite{CarmonaDelarue_book_II} for details.

As far as we know,  the earliest instance of such a generalization appeared in \cite{Huang} which proposed a linear-quadratic infinite-horizon model with a major player. Soon after, the finite-horizon counterpart of the model was considered in \cite{NguyenHuang1} and a first generalization to nonlinear cases was proposed in \cite{NourianCaines}. We believe theses are the first models of what is now called '\emph{mean field games with major and minor players'}. Still, the state of the major player does not enter the dynamics of the minor players, it only appears in their cost functionals. Later on \cite{NguyenHuang2} discussed a new approach to linear quadratic games  in which the major player's state enters the dynamics of the minor players. The authors solve the limiting control problem for the major player using a trick they call ``anticipative variational calculation''. 

The asymmetry between major and minor players was emphasized in \cite{BensoussanChauYam}
where the authors insist on the fact that the statistical distribution of the state  of a generic minor player should be derived endogenously.  Like in \cite{NourianCainesMalhame}, 
the paper \cite{BensoussanChauYam} characterizes the limiting problem by a set of stochastic partial differential equations. While working with the open loop formulation of the problem, the more recent account \cite{CarmonaZhu} also insists on the endogenous nature of the statical distribution of the state of a generic minor player. In fact, it goes one step further by reformulating the Mean Field Game with major and minor players as the search for a Nash equilibrium in a two player game 
over the time evolutions of states, some of which being of a McKean-Vlasov type. Note that, despite the fact that they offer a formal discussion of the general case, both papers  \cite{BensoussanChauYam} and \cite{CarmonaZhu} can only provide solutions in the linear quadratic case. For the sake of completeness, we also mention the recent technical report 
\cite{JaimungalNourian} where a major player is added to a particular case of the extended (in the sense that the interaction is through the controls) mean field game model of optimal execution introduced in Chapter 1 and solved in  Chapter 4 of \cite{CarmonaDelarue_book_I}.
Because of the absence of idiosyncratic noise, the initial conditions of the minor player states are assumed to be independent identically distributed random variables. The authors formulate a fixed point equilibrium problem when the rate of trading of the major player is restricted to be a linear function of the average rate of trading of the minor players, and they solve this fixed point problem
with deterministic controls in the infinite horizon stationary case.

In this paper, we present an alternative formulation for the Mean Field Games with major and minor players. In this approach, the search for Nash equilibria is naturally framed as the search for fixed points for the best response function for both types of players. 
As a \emph{fringe benefit} we are able to formulate and tackle the open and closed loop versions of the problem in one go.  
Beyond the fact that \cite{BensoussanChauYam} seems to be dealing only with the closed loop formulation of the problem, 
the main difference is the fact that instead of looking for a global Nash equilibrium of the whole system, including major and minor players, the authors choose a Stackelberg game strategy in which the major player goes first and  chooses its own control to minimize its expected cost, assuming that the response of the minor players to the choice of its control will be to put themselves in the (hopefully unique) mean field game equilibrium in the random environment induced by the control of the major player.
As a result, the finite-player game which is actually solved in \cite{BensoussanChauYam}, is merely a $N$-player game including only the minor players. In particular, the associated propagation of chaos is just a randomized version of the usual propagation of chaos associated to the usual mean field games. Here we follow the same line of attack as in \cite{CarmonaZhu}, making sure that the approximate equilibria obtained for finite player games are in fact $(N+1)$-player game equilibria including the major player as well as the $N$ minor players. 

\vskip 4pt
The paper is structured as follows. Our formulation of mean field games with major and minor players is presented in Section \ref{se:alternative}
below. There, we emphasize that as it relies on a fixed point argument in spaces of controls, and we explain how this approach can be used to tackle all sorts of versions of the game, whether the search is for open or closed loop (or even Markovian) equilibria.  Next, Section \ref{se:lq} implements this approach in the case of linear quadratic models. We recover the open loop solution of \cite{CarmonaZhu}, and provide a solution for closed loop models.
Section \ref{se:application} concludes with the solution of a generalization including a major player to the mean field game formulation proposed in \cite{NourianCainesMalhame} of a flocking model originally credited to Cucker and Smale \cite{CuckerSmale}. There, the dynamics of a large population of agents are governed by forces depicting the mean reversion of individual velocity to the mean velocity of the population. 
While early models of flocking do not involve any form of central coordination, several authors recently propose generalization of the flocking model by introducing leaders in the population. Such leaders have a pivotal impact on the rest of the population. In this spirit, we extend the mean field game formulation of \cite{NourianCainesMalhame} to include a major player which in equilibrium, should act as a free-will leader. We solve this model in the linear quadratic case, and we provide numerical simulations of the solution.

\section{Alternative Formulations for Mean Field Games with Major and Minor Players}
\label{se:alternative}

The goal of this section is to formulate the search for Nash equilibria for mean field games with major and minor players as a fixed point problem on a space of admissible controls. Since our discussion remains at the formal level, we do not introduce these mean field game models as limits of finite player games. We shall do just that only in the case of the linear quadratic models which we solve explicitly in Section \ref{se:lq} below. For pedagogical reasons, we treat separately the open and closed loop problems.
The rationale for this decision comes from the fact that, while solutions to the open and closed loop versions of the standard games often coincide in the mean field limit, this does not seem to be the case for games with major and minor players. Indeed, the characteristics of the state of the major player do not disappear in the limit when the number of minor players tends to infinity.
We shall illustrate this fact in our discussion of the linear quadratic models below.

\vskip 4pt
The general set up of a mean field game with major and minor players is as follows. The dynamics of the state of the system are given by stochastic differential equations of the form:
\begin{equation}
\label{fo:mmmfg_dyn}
\begin{cases}
dX^0_t&=b_0(t,X^0_t,\mu_t,\alpha^0_t)dt +\sigma_0(t,X^0_t,\mu_t,\alpha^0_t) dW^0_t\\
dX_t&=b(t,X_t,\mu_t,X^0_t,\alpha_t,\alpha^0_t)dt +\sigma(t,X_t,\mu_t,X^0_t,\alpha_t,\alpha^0_t dW_t,
\end{cases}
\end{equation}
where $\bW^0=(W^0_t)_{0\le t\le T}$ and $\bW=(W_t)_{0\le t\le T}$ are independent Wiener processes in $\RR^{d_0}$ and $\RR^d$ respectively, the quantities $X^0_t$, $\alpha^0_t$ with a superscript $0$ representing the state and the control of the major player while the he quantities $X_t$, $\alpha_t$ without a superscript represent the state and the control of the representative minor player. 
The controls $\alpha^0_t$ and $\alpha_t$ take values in closed convex subsets $A_0$ and $A$ of Euclidean spaces $\RR^{k_0}$ and $\RR^k$.
Here $\bmu=(\mu_t)_{0\le t\le T}$ is a measure valued process which in equilibrium, is expected to be given by the conditional distributions of the state of the representative minor player given the filtration $\FF^0=(\cF^0_t)_{0\le t\le T}$ generated by the Wiener process $\bW^0$ driving the dynamics of the state of the major player. Indeed, $\mu_t$ should be understood as a proxy for the empirical measure $\overline\mu^N_t$ of the states of $N$ minor players in the limit $N
\to\infty$. This limit is expected to be $\mu_t=\PP_{X_t|W^0_{[0,t]}}=\cL(X_t|W^0_{[0,t]})$
the conditional distribution of the state of the representative minor player given the initial path $W^0_{[0,t]}$ of the noise common to all the minor players,
namely the noise term driving the equation for the state of the major player. For later reference, we shall denote by $\FF=(\cF_t)_{0\le t\le T}$ the filtration generated by both Wiener processes.

The costs the players try to minimize are of the form:
\begin{equation}
\label{fo:mmmfg_costs}
\begin{cases}
J^0(\balpha^0,\balpha)&=\EE\bigl[\int_0^Tf_0(t,X^0_t,\mu_t,\alpha^0_t)dt +g^0(X^0_T,\mu_T)\bigr]\\
J(\balpha^0,\balpha)&=\EE\bigl[\int_0^Tf(t,X_t,\mu^N_t,X^0_t,\alpha_t,\alpha^0_t)dt +g(X_T,\mu_T)\bigr],
\end{cases}
\end{equation}
for some running and terminal cost functions $f_0$, $f$, $g_0$ and $g$.
The crucial feature of mean field games with major and minor players is that the dynamics of the state and the costs of the major player depend upon the statistical distribution of the states of the minor players while the states and the costs of the minor players depend upon not only their own states and the statistical distribution of the states of all the minor players, but also on the state and the control of the major player. This is what makes the analysis of these games more difficult than the standard mean field game models.

\vskip 4pt
We first treat the case of open loop equilibria for which we take advantage of the fact that the filtrations are assumed to be generated by the Wiener processes, to write the controls as functions of the paths of these Wiener processes. 

\subsubsection*{\textbf{Open Loop Version of the MFG Problem}}
\vskip -6pt

Here, we assume that the controls used by the major player and the representative minor player are of the form:
\begin{equation}
\label{fo:mmmfg_controls} 
\alpha^0_t=\phi^0(t,W^0_{[0,T]}),\quad\text{and}\quad  \alpha_t=\phi(t,W^0_{[0,T]},W_{[0,T]}),
\end{equation}
for deterministic progressively measurable functions 
$\phi^0:[0,T]\times\cC([0,T];\RR^{d_0})\mapsto A_0$ and 
$\phi:[0,T]\times\cC([0,T];\RR^{d})\times\cC([0,T];\RR^{d})\mapsto A$.
Progressive measurability of the function $\phi$ means that for each $t\in[0,T]$,  and $w^0,w\in\cC([0,T];\RR^{d})$, the value of 
$\phi(t,w^0,w)$ depends only upon the restrictions $w^0_{[0,t]}$ and $w_{[0,t]}$ of $w^0$ and $w$ to the interval $[0,t]$. Similarly for $\phi^0$. Our choice for the admissibility of the controls is consistent with our earlier discussion since we assume that the filtration $\FF^0$ and $\FF$ are generated by the Wiener processes $\bW^0$ and $(\bW^0,\bW)$ respectively.

\vskip 2pt
We understand a Nash equilibrium as a fixed point of the best response map. In the present context, the latter comprises two specific components: the best response of the major player to the behavior of all the minor players, and the best response of a representative minor player to the behavior of the major player and all the other minor players. So we need two separate steps to identify the best response map before we can define a Nash equilibrium as a fixed point of this map.

\vskip 6pt\noindent
\emph{The Major Player Best Response. }
We assume that the representative minor player uses the open loop control given by the progressively measurable function $\phi: (t,w^0,w)\mapsto \phi(t,w^0,w)$, so the problem of the major player is to minimize its expected cost:
\begin{equation}
\label{fo:mmmfg_major_cost}
J^{\phi,0}(\balpha^0)=\EE\Bigl[\int_0^Tf_0(t,X^0_t,\mu_t,\alpha^0_t)dt +g^0(X^0_T,\mu_T)\Bigr]
\end{equation}
under the dynamical constraints:
\begin{equation*}
\begin{cases}
dX^0_t&=b_0(t,X^0_t,\mu_t,\alpha^0_t)dt +\sigma_0(t,X^0_t,\mu_t,\alpha^0_t) dW^0_t\\
dX_t&=b(t,X_t,\mu_t,X^0_t,\phi(t,W^0_{[0,T]}, W_{[0,T]}),\alpha^0_t)dt +\sigma(t,X_t,\mu_t,X^0_t,\phi(t,W^0_{[0,T]}, W_{[0,T]}),\alpha^0_t) dW_t,
\end{cases}
\end{equation*}
where $\mu_t=\cL(X_t|W^0_{[0,t]})$ denotes the conditional distribution of $X_t$ given $W^0_{[0,t]}$. Since we are considering the open loop version of the problem, we search for minima in the class of controls $\balpha^0$ of the form $\alpha^0_t=\phi^0(t,W^0_{[0,T]})$ for a progressively measurable function $\phi^0$. So we frame the major player problem as the search for:
\begin{equation}
\label{fo:OL_MPP}
\phi^{0,*}(\phi)=\text{arg} \inf_{\balpha^0\leftrightarrow \phi^0} J^{\phi,0}(\balpha^0)
\end{equation}
where $\balpha^0\leftrightarrow \phi^0$ means that the infimum is over the set of controls $\balpha^0$ given by progressively measurable functions $\phi^0$.
For the sake of the present discussion, we assume implicitly that the argument of the minimization is not empty and reduces to a singleton. 
The important feature of this formulation is that the optimization of the major player appears naturally as an optimal control of the McKean-Vlasov type! In fact, it is an optimal control of the \emph{conditional} McKean-Vlasov type since the distribution appearing in the controlled dynamics is the conditional distribution of the state of the representative minor player.

\vskip 6pt\noindent
\emph{The Representative Minor Player Best Response. } 
To formulate the optimization problem of the representative minor player, we 
first describe the state of a system comprising a major player and a field of minor players different from the representative minor player we are focusing on.
So we assume that the major player uses a strategy $\balpha^0$ given by a progressively measurable function $\phi^0$ as in $\alpha^0_t=\phi^0(t,W^0_{[0,T]})$, and that the representative of the field of minor players uses a strategy $\balpha$ given by a progressively measurable function $\phi$ in the form $\alpha_t=\phi(t,W^0_{[0,T]},W_{[0,T]})$.
So the dynamics of the state of the system are given by:
\begin{equation*}
\begin{cases}
&dX^0_t=b_0(t,X^0_t,\mu_t,\phi^0(t,W^0_{[0,T]}))dt +\sigma_0(t,X^0_t,\mu_t,\phi^0(t,W^0_{[0,T]})) dW^0_t\\
&dX_t=b(t,X_t,\mu_t,X^0_t,\phi(t,W^0_{[0,T]},W_{[0,T]}),\phi^0(t,W^0_{[0,T]}))dt\\
&\hskip 125pt
 +\sigma(t,X_t,\mu_t,X^0_t,\phi(t,W^0_{[0,T]},W_{[0,T]}),\phi^0(t,W^0_{[0,T]})) dW_t,
\end{cases}
\end{equation*}
where as before, $\mu_t=\cL(X_t|W^0_{[0,t]})$ is the conditional distribution of $X_t$ given $W^0_{[0,t]}$. Notice that in the present situation, given the feedback functions $\phi^0$ and $\phi$, this stochastic differential equation in $\RR^{d_0}\times\RR^d$ giving the dynamics of the state of the system is of (conditional) McKean-Vlasov type since $\mu_t$ is the (conditional) distribution of (part of) the state.

\vskip 2pt
As explained earlier, we frame the problem of the representative minor player as the search for the best response to the major player and the field of the (other) minor players. So naturally, we formulate this best response as the result of the optimization problem of a virtual (extra) minor player which chooses a strategy $\overline{\balpha}$ given by a progressively measurable function $\overline\phi$ in the form $\overline\alpha_t=\overline\phi(t,W^0_{[0,T]},W_{[0,T]})$  in order to minimize its expected cost:
$$
J^{\phi^0,\phi}(\bar\balpha)=\EE\bigl[\int_0^T f(t,\overline X_t,X^0_t,\mu_t,\bar\alpha_t,\phi^0(t,W^0_{[0,T]}))dt +g(\overline X_T,\mu_t)\bigr],
$$
where the dynamics of the virtual state $\overline X_t$ are given by:
\begin{equation*}
\begin{split}
&d\overline X_t=b(t,\overline X_t,\mu_t,X^0_t,\bar\phi(t,W^0_{[0,T]},W_{[0,T]}),\phi^0(t,W^0_{[0,T]}))dt\\
&\hskip 125pt
 +\sigma(t,\overline X_t,\mu_t,X^0_t,\bar\phi(t,W^0_{[0,T]},W_{[0,T]}),\phi^0(t,W^0_{[0,T]})) d\overline W_t,
\end{split}
\end{equation*}
for a Wiener process $\overline{\bW}=(\overline W_t)_{0\le t\le T}$ independent of the other Wiener processes.
Notice that this optimization problem is not of McKean-Vlasov type. It is merely a classical optimal control problem, though with random coefficients.
As stated above, we search for minima in the class of feedback controls $\overline\balpha$ of the form $\overline\alpha_t=\overline\phi(t,W^0_{[0,T]},W_{[0,T]})$. We denote by:
\begin{equation}
\label{fo:OL_mPP}
\overline\phi^{*}(\phi^0,\phi)=\text{arg} \inf_{\overline{\balpha}\leftrightarrow \overline\phi} J^{\phi^0,\phi}(\bar\balpha)
\end{equation}
the result of the optimization.
Again, we assume that the optimal control exists, is given by a progressively measurable function, and is unique for the sake of convenience.

\vskip 4pt
We now formulate the existence of a Nash equilibrium for the mean field game with major and minor player as a fixed point of the best response maps
identified above by its components \eqref{fo:OL_MPP} and \eqref{fo:OL_mPP}. 
So by definition, a couple $(\hat\balpha^{0}, \hat\balpha)$ of controls given by progressively measurable functions  $(\hat\phi^{0}, \hat\phi)$
as above is a Nash equilibrium for the mean field game with major and minor players if it satisfies the fixed point equation:

\begin{equation}
\label{fo:mmmfg_fixed_point}
(\hat\phi^0,\hat\phi)=\big(\phi^{0,*}(\hat\phi),\bar\phi^{*}(\hat\phi^0,\hat\phi)\big).
\end{equation}

\subsubsection*{\textbf{Closed Loop Version of the MFG Problem}}
\vskip -6pt

The way we rewrote the open loop version of the problem may have been rather pompous, but it makes it easy to introduce the closed loop and Markovian versions of the problem. 
In this subsection, we assume that the controls used by the major player and the representative minor player are of the form:
$$ 
\alpha^0_t=\phi^0(t,X^0_{[0,T]},\mu_t),\quad\text{and}\quad  \alpha_t=\phi(t,X_{[0,T]},\mu_t,X^0_{[0,T]}), \quad  i=1,\cdots,N.
$$
for deterministic progressively measurable functions 
$\phi^0:[0,T]\times\cC([0,T];\RR^{d_0})\mapsto A_0$ and 
$\phi:[0,T]\times\cC([0,T];\RR^{d})\times\cC([0,T];\RR^{d})\mapsto A$.
The state $X^0_t$ of the major player and the state $X_t$ of the representative minor player evolve according to the same dynamic equations
\eqref{fo:mmmfg_dyn} as before, and the costs are also given by the same formula \eqref{fo:mmmfg_costs}, with $\mu_t=\cL(X_t|W^0_{[0,t]})$. We follow the same strategy as above to define the closed loop Nash equilibria of the game.

\vskip 6pt\noindent
\emph{The Major Player Best Response. }
We assume that the representative minor player uses the progressively measurable feedback function $\phi: (t,x,\mu,x^0)\mapsto \phi(t,x,\mu,x^0)$, so the problem of the major player is to minimize its expected cost \eqref{fo:mmmfg_major_cost}
under the dynamical constraints:
\begin{equation*}
\begin{cases}
dX^0_t&=b_0(t,X^0_t,\mu_t,\alpha^0_t)dt +\sigma_0(t,X^0_t,\mu_t,\alpha^0_t) dW^0_t\\
dX_t&=b(t,X_t,\mu_t,X^0_t,\phi(t,X_{[0,T]},\mu_t,X_{[0,T]}^0),\alpha^0_t)dt +\sigma(t,X_t,\mu_t,X^0_t,\phi(t,X_{[0,T]},\mu_t,X_{[0,T]}^0),\alpha^0_t) dW^i_t,
\end{cases}
\end{equation*}
whereas before  $\mu_t=\cL(X_t|W^0_{[0,t]})$ denotes the conditional distribution of $X_t$ given $W^0_{[0,t]}$. As explained earlier, we search for minima in the class of feedback controls $\balpha^0$ of the form $\alpha^0_t=\phi^0(t,X^0_{[0,T]},\mu_t)$, so we frame the major player problem as:
\begin{equation}
\label{fo:CL_MPP}
\phi^{0,*}(\phi)=\text{arg} \inf_{\balpha^0\leftrightarrow \phi^0} J^{\phi,0}(\balpha^0)
\end{equation}
which is an optimal control of the conditional McKean-Vlasov type! 

\vskip 6pt\noindent
\emph{The Representative Minor Player Best Response. } 
To formulate the optimization problem of the representative minor player, we 
first describe a system to which it needs to respond optimally. So we
assume that the major player uses a strategy $\balpha^0$ in feedback form given by a feedback function 
$\phi^0$ so that $\alpha^0_t=\phi^0(t,X^0_{[0,T]}\mu_t)$, and that the representative of the field of minor players uses a strategy $\balpha$ given by a 
progressively measurable feedback function $\phi$ in the form $\alpha_t=\phi(t,X_{[0,T]},X_{[0,T]}^0,\mu_t)$.
So the dynamics of the state of this system are given by:
\begin{equation*}
\begin{cases}
&dX^0_t=b_0(t,X^0_t,\mu_t,\phi^0(t,X^0_{[0,T]},\mu_t))dt +\sigma_0(t,X^0_t,\mu_t,\phi^0(t,X^0_{[0,T]},\mu_t)) dW^0_t\\
&dX_t=b(t,X_t,\mu_t,X^0_t,\phi(t,X_{[0,T]},X_{[0,T]}^0,\mu_t),\phi^0(t,X^0_{[0,T]},\mu_t))dt \\
&\hskip 125pt
+\sigma(t,X_t,\mu_t,X^0_t,\phi(t,X_{[0,T]},X_{[0,T]}^0,\mu_t),\phi^0(t,X^0_{[0,T]},\mu_t)) dW_t,
\end{cases}
\end{equation*}
where as before, $\mu_t=\cL(X_t|W^0_{[0,t]})$ is the conditional distribution of $X_t$ given $W^0_{[0,t]}$. Again, given the feedback functions $\phi^0$ and $\phi$, this stochastic differential equation in $\RR^{d_0}\times\RR^d$ is of (conditional) McKean-Vlasov type.

\vskip 2pt
As expected, we formulate this best response of the representative minor player as the result of the optimization problem of a virtual (extra) minor player which chooses a strategy $\overline{\balpha}$ given by a feedback function $\overline\phi$ in the form $\overline\alpha_t=\overline\phi(t,\overline X_t,X^0_t,\mu_t)$  in order to minimize its expected cost:
$$
J^{\phi^0,\phi}(\bar\balpha)=\EE\Bigl[\int_0^T f(t,\overline X_t,X^0_t,\mu_t,\bar\alpha_t,\phi^0(t,X^0_{[0,T]}\mu_t))dt +g(\overline X_T,\mu_t)\Bigr],
$$
where the dynamics of the virtual state $\overline X_t$ are given by:
\begin{equation*}
\begin{split}
&d\overline X_t=b(t,\overline X_t,\mu_t,X^0_t,\bar\phi(t,\overline X_{[0,T]},X_{[0,T]}^0,\mu_t),\phi^0(t,X^0_{[0,T]},\mu_t))dt \\
&\hskip 125pt
+\sigma(t,\overline X_t,\mu_t,X^0_t,\bar\phi(t,\overline X_{[0,T]},X_{[0,T]}^0,\mu_t),\phi^0(t,X^0_{[0,T]},\mu_t)) d\overline W_t,
\end{split}
\end{equation*}
for a Wiener process $\overline{\bW}=(\overline W_t)_{0\le t\le T}$ independent of the other Wiener processes.
We search for minima in the class of feedback controls $\overline\balpha$ of the form $\overline\alpha_t=\overline\phi(t,\overline X_{[0,T]},\mu_t,X^0_{[0,T]})$, and we denote the solution by:
\begin{equation}
\label{fo:CL_mPP}
\overline\phi^{*}(\phi^0,\phi)=\text{arg} \inf_{\overline{\balpha}\leftrightarrow \overline\phi} J^{\phi^0,\phi}(\bar\balpha).
\end{equation}
Since the best response map is given by its components \eqref{fo:CL_MPP} and \eqref{fo:CL_mPP}, we define the solution of a Nash equilibrium for the closed loop mean field game with major and minor player as the solution of the same fixed point equation \eqref{fo:mmmfg_fixed_point}, except for the fact that the functions $(\hat\phi^{0}, \hat\phi)$ are now progressively measurable
feedback functions of the type considered here.

\subsubsection{\textbf{Markovian Version of the MFG Problem}}
\vskip -6pt

Here, we assume that the controls used by the major player and the representative minor player are of the form:
$$ 
\alpha^0_t=\phi^0(t,X^0_t,\mu_t),\quad\text{and}\quad  \alpha_t=\phi(t,X_t,\mu_t,X^0_t), \quad  i=1,\cdots,N.
$$
for deterministic feedback functions $\phi^0:[0,T]\times\RR^{d_0}\times\cP_2(\RR^d)\mapsto A_0$ and $\phi:[0,T]\times\RR^{d}\times\cP_2(\RR^d)\times\RR^{d_0}\mapsto A$.
The state $X^0_t$ of the major player and the state $X_t$ of the representative minor player evolve according to the same dynamic equations
\eqref{fo:mmmfg_dyn} as before and the costs are also given by the same formula \eqref{fo:mmmfg_costs}, with $\mu_t=\cL(X_t|W^0_{[0,t]})$.

\vskip 6pt\noindent
\emph{The Major Player Best Response. }
We assume that the representative minor player uses the feedback function $\phi: (t,x,\mu,x^0)\mapsto \phi(t,x,\mu,x^0)$, so the problem of the major player is to minimize its expected cost \eqref{fo:mmmfg_major_cost}
under the dynamical constraints:
\begin{equation*}
\begin{cases}
dX^0_t&=b_0(t,X^0_t,\mu_t,\alpha^0_t)dt +\sigma_0(t,X^0_t,\mu_t,\alpha^0_t) dW^0_t\\
dX_t&=b(t,X_t,\mu_t,X^0_t,\phi(t,X_t,\mu_t,X_t^0),\alpha^0_t)dt +\sigma(t,X_t,\mu_t,X^0_t,\phi(t,X_t,\mu_t,X_t^0),\alpha^0_t) dW^i_t,
\end{cases}
\end{equation*}
where as before $\mu_t=\cL(X_t|W^0_{[0,t]})$ denotes the conditional distribution of $X_t$ given $W^0_{[0,t]}$. We search for minima in the class of feedback controls $\balpha^0$ of the form $\alpha^0_t=\phi^0(t,X^0_t,\mu_t)$, so we frame the major player problem as:
\begin{equation}
\label{fo:M_MPP}
\phi^{0,*}(\phi)=\text{arg} \inf_{\balpha^0\leftrightarrow \phi^0} J^{\phi,0}(\balpha^0)
\end{equation}
As before, the optimization problem of the major player is of the conditional Mckean-Vlasov type.

\vskip 6pt\noindent
\emph{The Representative Minor Player Best Response. } 
To formulate the optimization problem of the representative minor player, we 
first describe a system to which it needs to respond optimally. So we
assume that the major player uses a strategy $\balpha^0$ in feedback form given by a feedback function 
$\phi^0$ so that $\alpha^0_t=\phi^0(t,X^0_t,\mu_t)$, and that the representative of the field of minor players uses a strategy $\balpha$ given by a feedback function $\phi$ in the form $\alpha_t=\phi(t,X_t,X_t^0,\mu_t)$.
So the dynamics of the state of this system are given by:
\begin{equation*}
\begin{cases}
&dX^0_t=b_0(t,X^0_t,\mu_t,\phi^0(t,X^0_t,\mu_t))dt +\sigma_0(t,X^0_t,\mu_t,\phi^0(t,X^0_t,\mu_t)) dW^0_t\\
&dX_t=b(t,X_t,\mu_t,X^0_t,\phi(t,X_t,X_t^0,\mu_t),\phi^0(t,X^0_t,\mu_t))dt +\sigma(t,X_t,\mu_t,X^0_t,\phi(t,X_t,X_t^0,\mu_t),\phi^0(t,X^0_t,\mu_t)) dW_t,
\end{cases}
\end{equation*}
where as before, $\mu_t=\cL(X_t|W^0_{[0,t]})$ is the conditional distribution of $X_t$ given $W^0_{[0,t]}$. Again, given the feedback functions $\phi^0$ and $\phi$, this stochastic differential equation in $\RR^{d_0}\times\RR^d$ is of (conditional) McKean-Vlasov type.

\vskip 2pt
As before, we frame the problem of the representative minor player as the search for the best response to the behavior of the major player and the field of the (other) minor players. So we solve the optimization problem of a virtual (extra) minor player which chooses a strategy $\overline{\balpha}$ given by a feedback function $\overline\phi$ in the form $\overline\alpha_t=\overline\phi(t,\overline X_t,X^0_t,\mu_t)$  in order to minimize its expected cost:
$$
J^{\phi^0,\phi}(\bar\balpha)=\EE\Bigl[\int_0^T f(t,\overline X_t,X^0_t,\mu_t,\bar\alpha_t,\phi^0(t,X^0_t,\mu_t))dt +g(\overline X_T,\mu_t)\Bigr],
$$
where the dynamics of the virtual state $\overline X_t$ are given by:
\begin{equation*}
\begin{split}
&d\overline X_t=b(t,\overline X_t,\mu_t,X^0_t,\bar\phi(t,\overline X_t,X_t^0,\mu_t),\phi^0(t,X^0_t,\mu_t))dt\\
&\hskip 125pt
 +\sigma(t,\overline X_t,\mu_t,X^0_t,\bar\phi(t,\overline X_t,X_t^0,\mu_t),\phi^0(t,X^0_t,\mu_t)) d\overline W_t,
\end{split}
\end{equation*}
for a Wiener process $\overline{\bW}=(\overline W_t)_{0\le t\le T}$ independent of the other Wiener processes.
We search for minima in the class of feedback controls $\overline\balpha$ of the form $\overline\alpha_t=\overline\phi(t,\overline X_t,\mu_t,X^0_t)$, and we denote the solution by:
\begin{equation}
\label{fo:M_mPP}
\overline\phi^{*}(\phi^0,\phi)=\text{arg} \inf_{\overline{\balpha}\leftrightarrow \overline\phi} J^{\phi^0,\phi}(\bar\balpha).
\end{equation}
Finally, we define the solution of a Nash equilibrium for the Markovian mean field game with major and minor player as the solution of the same fixed point equation \eqref{fo:mmmfg_fixed_point}, except for the fact that the functions $(\hat\phi^{0}, \hat\phi)$ are now
feedback functions of the type considered here.

\section{Linear Quadratic Models}
\label{se:lq}

In this section, we consider the mean field game with major and minor players issued from the finite player game in which the dynamics of the states of the players are given by the following linear stochastic differential equations:
\begin{equation}
\label{fo:finite_players}
\begin{dcases}
dX^{N,0}_t=({L}_0 X^{N,0}_t+B_0 \alpha^{N,0}_t+F_0 \bar{X}^N_t)dt+D_0 dW^0_t,
\\
dX^{N,i}_t=(L X^{N,i}_t+B\alpha^{N,i}_t+F\bar{X}^N_t+GX^0_t)dt+DdW^i_t,\qquad 1\le i\le N,
\end{dcases}
\end{equation}
for $t \in [0,T]$, and we choose $A_{0}=\RR^{k_{0}}$
and $A=\RR^k$.  
The coefficients are deterministic constant matrices independent of time. 
The real matrices $L_0$, $B_0$, $F_0$ and $D_0$ are of dimensions $d_0\times d_0$, $d_0\times k_0$, $d_0\times d$ and $d_0\times m_0$ respectively. Similarly, the real matrices $L$, $B$, $F$, $G$ and $D$ are of dimensions $d\times d$, $d\times k$, $d\times d$, $d\times d_0$, and $d_0\times m_0$ respectively.
The cost functionals for the major and minor players are given by:
\begin{equation*}
\begin{split}
&J^{N,0}\bigl(\balpha^{N,0},\cdots,\balpha^{N,N}\bigr)
\\
&\hspace{15pt}=\EE\bigg[\int^T_0
\Big[\bigl(X^{N,0}_t-\Psi_0(\bar{X}^N_t)\bigr)^\dagger Q_0 
\bigl(X^{N,0}_t-\Psi_0(\bar{X}^N_t)\bigr)+(\alpha^{N,0}_{t})^\dagger R_0 \alpha^{N,0}_t\Big]\,dt\bigg],
\\
&J^{N,i}\bigl(\balpha^{N,0},\cdots,\balpha^{N,N}\bigr)
\\
&\hspace{-5pt}=\mathbb{E}\bigg[\int^T_0 \Big[\bigl(X^{N,i}_t-\Psi(X^{N,0}_t,\bar{X}^N_t)\bigr)^\dagger Q 
\bigl(X^{N,i}_t-\Psi(X^{N,0}_t,\bar{X}^N_t))+(\alpha^{N,i}_t)^\dagger R \alpha^{N,i}_t\Big]dt\bigg],
\end{split}
\end{equation*}
in which $Q_0$, $Q$, $R_0$ and $R$ are positive definite symmetric matrices of dimensions $d_0\times d_0$, $d\times d$, $k_0\times k_0$ and $k\times k$, and where the functions $\Psi_0$ and $\Psi$ are defined by:
$$
\Psi_0(X)=H_0 X+\eta_0,\qquad \Psi(X,Y)=H X+H_1 Y+\eta,
$$
for some fixed $d_0\times d$, $d\times d_0$ and $d\times d$ matrices $H_0$, $H$ and $H_1$, and some fixed $\eta_0\in\RR^{d_0}$ and $\eta\in\RR^d$. Here, $\bar{X}^N_t$ stands for the empirical mean $(X^{N,1}_t+\cdots + X^{N,N}_t)/N$.

\vskip 4pt
We chose to study this specific linear quadratic model to match existing literature on the subject. Several variants are possible which can be treated using the same procedure. See for example the application discussed in Section \ref{se:application} below.

\subsubsection*{\textbf{Open-Loop Equilibrium}}
\vskip -6pt
In the mean field limit, the dynamics \eqref{fo:finite_players} of the major player state $X_t^0$ and the state $X_t$ of the representative minor player are given by:
\begin{equation}
\label{fo:mmmfg_alternative_ol_state}
\begin{dcases}
dX_t^0\;\;&= (L_0 X_t^0 + B_0 \alpha_t^0 + F_0 \bar X_t) dt + D_0 dW_t^0\\
dX_t\;\; &= (L X_t + B \alpha_t + F \bar X_t + G X_t^0) dt + D dW_t
\end{dcases}
\end{equation}
where $\bar X_t = \mathbb{E}[X_t | \mathcal{F}_t^0]$ is the conditional expectation of $X_t$ with respect to the filtration generated by the history of the Wiener process $\bW^0$ up to time $t$. Accordingly, the cost functionals for the major and minor players are given by:
\begin{align*}
&J^0(\balpha^0,\balpha) = \EE\left [ \int_{0}^T [(X_t^0 - H_0 \bar X_t - \eta_0)^{\dagger} Q_0 (X_t^0 - H_0 \bar X_t - \eta_0) + \alpha_t^{0\dagger}R_0 \alpha_t^0] dt\right]\\
&J(\balpha^0,\balpha)= \EE\left [ \int_{0}^T [(X_t - H X_t^0 - H_1 \bar X_t -\eta)^{\dagger} Q (X_t - H X_t^0 - H_1 \bar X_t -\eta) + \alpha_t^{\dagger}R \alpha_t] dt\right]
\end{align*}
in which $Q$, $Q_0$, $R$, $R_0$ are symmetric matrices, and $R$, $R_0$ are assumed to be positive definite.
Taking conditional expectations in the equation for the state of the representative minor player we get:
\begin{equation}
\label{fo:mmmfg_cond_exp}
d\bar X_t= [(L + F) \bar X_t + B  \overline\alpha_t + G X_t^0]\;dt,
\end{equation}
with $\overline\alpha_t=\EE[\alpha_t|\cF_t^0]$. The idea is now to express the optimization problem of the major player over the dynamics of the couple $(\overline X_t,X^0_t)$. In order to do so, we introduce the following notation:
$$
\begin{array}{c}
\mathbb{X}_t = \left[\begin{array}{c}\bar X_t \\ X_t^0 \end{array}\right],\;\;\mathbb{L}_0 =  \left[\begin{array}{cc}L+F&G\\ F_0 &L_0 \end{array}\right],\;\;\mathbb{B}_0 = \left[\begin{array}{c}0 \\ B_0 \end{array}\right],\;\;\mathbb{B} = \left[\begin{array}{c}B \\ 0 \end{array}\right], \mathbb{D}_0 = \left[\begin{array}{c}0 \\ D_0 \end{array}\right]\\\\
\mathbb{F}_0 = \left[\begin{array}{cc}H_0^\dagger Q_0 H_0 & - H_0^\dagger Q_0 \\ -Q_0H_0 & Q_0 \end{array}\right],\;\;f_0= \left[\begin{array}{c}H_0^\dagger Q_0 \eta_0 \\ -Q_0\eta_0 \end{array}\right].
\end{array}
$$
Notice that, the fact that the matrix $Q_0$ is symmetric non-negative definite implies that $\FF_0$ is also symmetric non-negative definite.
This will play a crucial role when we face the solution of certain matrix Riccati equations. The optimization problem of the major player becomes:
$$
\inf_{\balpha^0\in\AA_0}\mathbb{E}\left[\int_{0}^T [ \mathbb{X}_t^\dagger \mathbb{F}_0\mathbb{X}_t +2 \mathbb{X}_t^\dagger f_0 +  \eta_0^\dagger Q_0 \eta_0 + \alpha_t^{0\dagger}R_0 \alpha_t^0] dt  \right],
$$
where the controlled dynamics are given by:
\begin{equation}
\label{fo:mmmfg_semi_state}
d\mathbb{X}_t = (\mathbb{L}_0 \mathbb{X}_t +\mathbb{B}_0 \alpha_t^0 +  \mathbb{B}  \overline\alpha_t) dt + \mathbb{D}_0 dW_t^0.
\end{equation}
The reduced Hamiltonian is given by:
\[
H^{(r),\overline\alpha}(t, x, y, \alpha^0) = y^\dagger (\mathbb{L}_0 x +\mathbb{B}_0 \alpha^0 +  \mathbb{B} \overline\alpha_t) + x^\dagger \mathbb{F}_0 x +2  x^\dagger f_0 +  \eta_0^\dagger Q_0 \eta_0 +\alpha^{0\dagger}R_0 \alpha^0.
\]
Here we added the superscript $\overline\alpha$ for the Hamiltonian in order to emphasize that the optimization of the major player is performed assuming that the representative minor player is using strategy $\balpha  \in \AA$. Obviously, $H^{(r),\overline\alpha}$ is a random function, the randomness coming from the realization of the control of the representative minor player. However we see that almost surely $\RR^{d_0+d}\times A_0\ni (x, \alpha^0)\rightarrow H^{(r),\overline\alpha}(t, x, y, \alpha^0)$ is jointly convex, and we can use the sufficient condition of the stochastic maximum principle. Therefore the minimizer of the reduced Hamiltonian and the optimal control are given by:  
$$
\hat\alpha^{0} = - \frac12 R_0^{-1} \mathbb{B}^{0\dagger} y,
\qquad\text{and}\qquad
\hat\alpha^0_t = - \frac12 R_0^{-1} \mathbb{B}^{0\dagger} \YY_t,
$$
respectively, where $(\XX_t, \YY_t)_{0\le t\le T}$ solves the forward-backward stochastic differential equation:
\begin{equation}
\label{fo:mmmfg_lq_fbsde_major}
\begin{dcases}
d\XX_t &=\;\; (\LL_0 \mathbb{X}_t -\frac12 \BB_0 R_0^{-1}\BB_{0}^{\dagger}\YY_t +  \BB \overline\alpha_t) dt + \DD_0 dW_t^0\\
d\YY_t &=\;\; -(\LL_{0}^{\dagger} \YY_t +2\FF_{0}\XX_t  + 2f_0) dt + \ZZ_t dW_t^0,\;\;\;\YY_T =0.
\end{dcases}
\end{equation}

\vskip 4pt
We now consider the representative minor player's problem. We fix an admissible strategy $\balpha^0\in\AA_0$ for the major player, and an admissible strategy $\balpha\in\AA$ for the representative of the \emph{other} minor players, and its $\FF^0$-optional projection $\overline\balpha$ defined by $\overline\alpha_t=\EE[\alpha_t|\cF^0_t]$. This prescription leads to the time evolution of the state of a system given by
\eqref{fo:mmmfg_alternative_ol_state}, equation \eqref{fo:mmmfg_cond_exp} after taking conditional expectations, and finally the dynamic equation
\eqref{fo:mmmfg_semi_state}. Given this background state evolution, 
the representative minor player needs to solve:
$$
\inf_{\tilde\balpha\in\AA}\EE \left[\int_0^T [(\tilde X_t - [H_1, H] \XX_t - \eta)^{\dagger} Q(\tilde X_t - [H_1, H] \XX_t - \eta)+\tilde\alpha_t^\dagger R \tilde\alpha_t ] dt \right],
$$
where the dynamics of the controlled state $\tilde X_t$ are given by:
$$
d\tilde X_t = (L \tilde X_t + B \tilde\alpha_t + [F,G]\XX_t) dt + D dW_t.
$$
Note that the process $\mathbb{X}_t$ is merely part of the random coefficients of the optimization problem. We introduce the reduced Hamiltonian:
\begin{equation*}
\begin{split}
&H^{(r),\alpha^0,\alpha}(t,\tilde x,\tilde y,\tilde\alpha) = 
\tilde y^{\dagger}(L \tilde x + B \tilde\alpha + [F,G] \mathbb{X}_t)\\ 
&\hskip 40pt
+ (\tilde x - [H_1, H] \mathbb{X}_t - \eta)^{\dagger} Q(\tilde x - [H_1,H] \mathbb{X}_t - \eta)+\tilde\alpha^\dagger R \tilde\alpha.
\end{split}
\end{equation*}
Once again we use the superscript $(\alpha^0,\alpha)$ to emphasize the fact that the optimization is performed under the environment created by the major player using strategy $\balpha^0$ and the population of the other minor players using $\balpha$, leading to the use of its $\FF^0$-optional projection $\overline\balpha$. $H^{(r),\alpha^0,\alpha}$ depends on the random realization of the environment and is almost surely jointly convex  in $(\tilde x,\tilde\alpha)$. Applying the stochastic maximum principle, the optimal control exists and is given by $\tilde\alpha_t = -\frac12 R^{-1}B^\dagger \tilde Y_t$, where $(\tilde{\bX},\tilde{\bY})$ solves the following FBSDE: 
\begin{equation}
\label{fo:mmmfg_lq_fbsde_minor}
\begin{dcases}
d\tilde X_t &=\;\; (L \tilde X_t -B R^{-1}B^{\dagger}\tilde Y_t +  [F,G]\tilde{\XX}_t) dt + D dW_t\\
d\tilde Y_t &=\;\; -\bigl(L^{\dagger} \tilde Y_t +2 Q\bigl( X_t -  [H_1, H] \tilde{\XX}_t - \eta\bigr)\bigr) dt + Z_t dW_t + Z_t^0 dW_t^0,
\end{dcases}
\end{equation}
with terminal condition $Y_T =0$.  Recall that in this FBSDE, the process $(\XX_t)_{0\le t\le T}$ only acts as a random coefficient.
It is determined \emph{off line} by solving the standard stochastic differential equation:
\begin{equation}
\label{fo:mmmfg_tilde}
d\tilde{\XX}_t =\;\; (\mathbb{L}_0 \tilde{\XX}_t +\mathbb{B}_0 \alpha_t^0 +  \mathbb{B} \overline\alpha_t) dt + \mathbb{D}_0 dW_t^0\\
\end{equation}
Notice that equation \eqref{fo:mmmfg_tilde} is exactly the same equation as \eqref{fo:mmmfg_semi_state}. Still, we use a different notation for the solution. Indeed, at this stage of the proof (i.e. before considering the fixed point step), the coefficient processes $(\alpha_t^0)_{0\le t\le T}$ and $(\overline\alpha_t)_{0\le t\le T}$ are (likely to be) different, preventing us from identifying the solutions of \eqref{fo:mmmfg_tilde} and \eqref{fo:mmmfg_semi_state}.

\vskip 4pt
Now that we are done characterizing the solutions of both optimization problems, we identify the fixed point constraint in the framework given by the characterizations of the two optimization problems, 
The fixed point condition \eqref{fo:mmmfg_fixed_point} characterizing Nash equilibria in the current set-up says that:
\begin{equation*}
\alpha_t^0 = - \frac12 R_0^{-1} \BB_{0}^{\dagger} \YY_t,
\end{equation*}
where $(\YY_t)_{0\le t\le T}$ is the backward component of the solution of \eqref{fo:mmmfg_lq_fbsde_major} with 
$\overline\alpha_t =\;\EE[\alpha_t|\mathcal{F}_t^0]$, and:
$$
\alpha_t=\tilde\alpha_t = -\frac12 R^{-1}B^\dagger \tilde Y_t,
$$ 
where $(\tilde{Y}_t)_{0\le t\le T}$ is the backward component of the solution of \eqref{fo:mmmfg_lq_fbsde_minor} in which the random coefficient
$(\tilde\XX_t)_{0\le t\le T}$ solves \eqref{fo:mmmfg_tilde} with the processes $(\alpha^0_t)_{0\le t\le T}$ and $(\overline\alpha_t)_{0\le t\le T}$ just defined.
So in equilibrium, equations \eqref{fo:mmmfg_tilde} and \eqref{fo:mmmfg_semi_state} have the same coefficients and we can identify their solutions 
$(\XX_t)_{0\le t\le T}$ and $(\tilde\XX_t)_{0\le t\le T}$.

\vskip 2pt
The optimal controls for the major and representative minor players are functions of the solution of the following FBSDE which we obtain by putting together the FBSDEs \eqref{fo:mmmfg_lq_fbsde_major} and \eqref{fo:mmmfg_lq_fbsde_minor} characterizing the major and representative minor players' optimization problem:
\begin{equation}
\label{fo:mmmfg_lq_fbsde_Nash}
\begin{dcases}
&d\XX_t = (\mathbb{L}_0 \XX_t -\frac12 \BB_0 R_0^{-1}\BB_{0}^{\dagger}\YY_t - \frac12\BB R^{-1}B^\dagger  \EE[\tilde Y_t|\cF_t^0]) dt + \DD_0 dW_t^0\\
&d\tilde X_t = (L \tilde X_t - \frac12 B R^{-1}B^{\dagger}\tilde Y_t +  [F,G]\XX_t) dt + D dW_t\\
&d\YY_t = -(\LL_{0}^{\dagger} \YY_t +\FF_{0}\XX_t  + f_0) dt + \ZZ_t dW_t^0,\;\;\;\YY_T =0\\
&d\tilde Y_t = -(L^{\dagger} \tilde Y_t + 2 Q\tilde X_t -2 Q [H_1, H] \XX_t - 2 Q\eta) dt + Z_t dW_t + Z_t^0 dW_t^0,\;\tilde Y_T =0.
\end{dcases}
\end{equation}
We summarize the above discussion in the form of a verification theorem for open-loop Nash equilibrium.

\begin{theorem}
If the system \eqref{fo:mmmfg_lq_fbsde_Nash} admits a solution, then the linear quadratic mean field game problem with major  and minor players admits an open-loop Nash equilibrium. The equilibrium strategy $(\balpha^0,\balpha)$ is given by $\hat\alpha_t^0 =  -(1/2)R_0^{-1} \mathbb{B}_{0}^{\dagger} \mathbb{Y}_t$ for the major player and $\hat\alpha_t = - (1/2) R^{-1}B \tilde Y_t$ for the representative minor player.
\end{theorem}

\vskip 4pt
The way the system \eqref{fo:mmmfg_lq_fbsde_Nash} is stated is a natural conclusion of the search for equilibrium as formulated by the fixed point
step following the two optimization problems. However, as convenient as can be, simple remarks can help the solution of this system.
First we notice one could solve for $(\XX_t,\YY_t)_{0\le t\le T}$ by solving the FBSDE formed by the first and the third equations if we knew 
$\overline Y_t =\EE[\tilde Y_t|\cF_t^0]$. By taking conditional expectations with respect to $\cF^0_t$ in the second equation, 
and by subtracting the result from the equation satisfied by the first component of the first equation, we identify $\EE[\tilde X_t|\cF_t^0]$ with $\overline X_t$ because they have the same initial conditions. Next, by taking conditional expectations with respect to $\cF^0_t$ in the fourth equation, we see that $(\overline Y_t)_{0\le t\le T}$ should satisfy:
$$
d\overline Y_t = -(L^{\dagger} \overline Y_t +Q\overline X_t - Q [H_1, H] \XX_t  - Q\eta) dt  + \overline Z_t^0 dW_t^0,\;\overline Y_T =0
$$
Consequently, the solution of \eqref{fo:mmmfg_lq_fbsde_Nash} also satisfies:
\begin{equation}
\label{fo:mmmfg_lq_fbsde_final}
\begin{dcases}
&d\XX_t = (\mathbb{L}_0 \XX_t -\frac12 \BB_0 R_0^{-1}\BB_{0}^{\dagger}\YY_t - \frac12\BB R^{-1}B^\dagger \overline Y_t) dt + \DD_0 dW_t^0\\
&d\YY_t = -(\LL_{0}^{\dagger} \YY_t +2\FF_{0}\XX_t  + 2f_0) dt + \ZZ_t dW_t^0,\;\;\;\YY_T =0\\
&d\overline Y_t = -\bigl(L^{\dagger} \overline Y_t + 2 \bigl([Q,0] - Q [H_1, H] \bigr)\XX_t  - 2 Q\eta\bigr) dt  + \overline Z_t^0 dW_t^0,\;\overline Y_T =0.
\end{dcases}
\end{equation}
Our final remark is that the solution of system \eqref{fo:mmmfg_lq_fbsde_final} is not only necessary, but also sufficient. Indeed, once it is solved, 
one can solve for $(\tilde X_t,\tilde Y_t)_{0\le t\le T}$ by solving the affine FBSDE with random coefficients formed by the second and fourth equations of \eqref{fo:mmmfg_lq_fbsde_Nash} and check that $\EE[\tilde Y_t|\cF^0_t]$ is indeed the solution of the third equation of \eqref{fo:mmmfg_lq_fbsde_final}.

\vskip 4pt
Identifying $\YY_t$ with $[\overline P_t^\dagger,\overline P_t^{0\dagger}]^\dagger$ we recognize the FBSDE used in \cite{CarmonaZhu}.

\subsubsection*{\textbf{A  Closed Loop Equilibrium}}
In this section we implement the closed loop alternative formulation of the equilibrium problem. Since we expect that the optimal controls will be in feedback form, we search directly for Markovian controls.
In other words, we assume that the controls used by major player and minor players are respectively of the form:
$$
\alpha_t^0 = \phi^0(t, X_t^0, \bar X_t),
\qquad\text{and}\qquad
\alpha_t = \phi(t, X_t, X_t^0, \bar X_t),
$$
for some $\RR^{k_0}$ and $\RR^k$ valued deterministic functions $\phi^0$ and $\phi$
defined on $[0,T]\times\RR^{d_0}\times \RR^d$ and $[0,T]\times\RR^d\times\RR^{d_0}\times \RR^d$ respectively.  For the sake of simplicity, we assume that $A_0= \RR^{k_0}$ and $A=\RR^k$. 
So the major player can only observe its own state and the mean of minor player's states, while the representative minor player can observe its own state, the state of the major player, as well as the mean of the other minor players' states. 
This version of the equilibrium problem is more difficult than its open loop analog. For that reason, we are not trying to 
construct the best response map for all the possible choices of control processes $\balpha^0$ and $\balpha$.
Instead, we construct it for a restricted class of feedback functions $\phi^0$ and $\phi$ in which we can still find a fixed point, hence a Nash equilibrium.

To be more specific, we construct the best responses to controls $\balpha^0$ and $\balpha$ of the form: 
\begin{align}
\alpha_t^0 &=\; \phi^0(t, X_t^0, \bar X_t) =  \phi^0_0(t) + \phi^0_1(t) X_t^0 + \phi^0_2(t) \bar X_t\label{fo:alpha0}\\
\alpha_t &=\;\;\phi(t, X_t, X_t^0, \bar X_t)=  \phi_0(t) + \phi_1(t) X_t + \phi_2(t) X_t^0 + \phi_3(t) \bar X_t\label{fo:alpha}
\end{align}
where the functions $[0,T]\ni t\rightarrow \phi^0_i(t)$ for  $i=0,1,2$ and $[0,T]\ni t \rightarrow \phi_i(t)$ for $i = 0,1,2,3$ are matrix-valued deterministic continuous functions with the appropriate dimensions, in other words, $\phi_0^0(t)\in\RR^{k_0}$, $\phi_1^0(t)\in\RR^{k_0\times d_0}$, $\phi_2^0(t)\in\RR^{k_0\times d}$, $\phi_0(t)\in\RR^{k}$, $\phi_1(t)\in\RR^{k\times d}$, $\phi_2(t)\in\RR^{k\times d_0}$, and $\phi_3(t)\in\RR^{k\times d}$. 

\vskip 4pt
We first consider the major player's optimization problem. We assume that the representative minor player uses strategy $\alpha_t=\phi(t, X_t, X_t^0, \bar X_t)$ as specified in \eqref{fo:alpha}. Next we look for the control $\balpha^0$ 
which could be used by the major player to minimize its expected cost. The dynamics of the system is then given by:
\begin{equation}
\label{fo:mmmfg_alternative_cl_state}
\begin{cases}
&\hskip -10pt
dX_t^0 = (L_0 X_t^0 + B_0 \alpha_t^0 + F_0 \bar X_t) dt + D_0 dW_t^0\\
&\hskip -10pt
dX_t =  \Bigl[ B\phi_0(t)  + (L+B\phi_1(t) ) X_t +  (B\phi_2(t)+G) X^0_t + (B\phi_3(t)+F) \overline X_t) \Bigr]dt + D dW_t,
\end{cases}
\end{equation}
where as before $\overline X_t = \mathbb{E}[X_t | \mathcal{F}_t^0]$ is the conditional expectation of $X_t$ with respect to the filtration generated by the history of the Wiener process $\bW^0$ up to time $t$.  In their current form, the dynamics of the couple $(X^0_t,X_t)$ are of a McKean-Vlasov type since the mean of $X_t$ appears in the coefficients of the equation giving $dX^0_t$. However, in
order to find a minimalist version of dynamical equations for a state over which the optimization problem of the major player can be formulated, we take conditional expectations in the equation for the state of the representative minor player. We get:
\begin{equation}
\label{fo:mmmfg_cl_cond_exp}
d\overline X_t=  \Bigl[ B\phi_0(t)  + (L+B[\phi_1(t) +\phi_3(t)] +F) \overline X_t +  (B\phi_2(t)+G) X^0_t \Bigr]dt.
\end{equation}
As in the case of the open loop version of the equilibrium problem, we express the optimization problem of the major player over the dynamics of the couple $(\overline X_t,X^0_t)$. In order to do so, we use the same notation $\XX_t$, $\FF_0$, $f_0$, $\BB_0$, $\BB$, $\DD$ and $\DD_0$ as in the case of our analysis of the open loop problem, and we introduce the following new notation:
$$
\LL^{(cl)}_0 (t)=  
\left[\begin{array}{cc}L+B[\phi_1(t)+\phi_3(t)]+F&B\phi_2(t)+G\\ F_0 &L_0 \end{array}\right],
\quad
\CC^{(cl)}_0 = \left[\begin{array}{c}B\phi_0(t) \\0 \end{array}\right],
$$
and the optimization problem of the major player can be formulated exactly as in the open loop case as the minimization:
$$
\inf_{\balpha^0\in\AA_0}\mathbb{E}\left[\int_{0}^T [ \mathbb{X}_t^\dagger \mathbb{F}_0\mathbb{X}_t +2 \mathbb{X}_t^\dagger f_0 +  \eta_0^\dagger Q_0 \eta_0 + \alpha_t^{0\dagger}R_0 \alpha_t^0] dt  \right]
$$
where the controlled dynamics are given by:
\begin{equation}
\label{fo:mmmfg_semi_state}
d\XX_t =\bigl[\LL^{(cl)}_0(t) \XX_t +\BB_0 \alpha_t^0 +  \CC^{(cl)}_0(t)\bigr] dt + \mathbb{D}_0 dW_t^0.
\end{equation}
The reduced Hamiltonian (minus the term $ \eta_0^\dagger Q_0 \eta_0$ which is irrelevant) is given by:
$$
H^{(r),\phi}(t, x, y, \alpha^0) = y^\dagger [\LL^{(cl)}_0 x +\BB_0 \alpha^0 +  \CC^{cl}(t)] + x^\dagger \mathbb{F}_0 x +2  x^\dagger f_0 + \alpha^{0\dagger}R_0 \alpha^0.
$$
Applying the stochastic maximum principle, we find that the optimal control is given as before by $\hat\alpha_t^{0} = - (1/2)R_0^{-1}\BB_0^\dagger \mathbb{Y}_t$, where $(\XX_t, \YY_t, \ZZ_t)_{0\le t\le T}$ solves the linear FBSDE:
\begin{equation}
\label{fo:fbsde_xx_yy}
\begin{cases}
d\XX_t &=\;[\LL_0^{(cl)}(t) \XX_t -\frac12 \BB_0 R_0^{-1}\BB_{0}^{\dagger}\YY_t +  \CC^{(cl)}_0(t)] dt + \mathbb{D}_0 dW_t^0\\
d\YY_t &=\; -[\LL_0^{(cl)}(t)^{\dagger} \YY_t +2\FF_{0}\mathbb{X}_t  + 2 f_0) dt + \mathbb{Z}_t dW_t^0,\;\;\;\mathbb{Y}_T =0.
\end{cases}
\end{equation}
This FBSDE being affine, we expect the decoupling field to be affine as well, so we search for a solution of the form $\YY_t = K_t\XX_t + k_t$ for two deterministic functions $t\mapsto K_t\in\RR^{(d+d_0)\times (d+d_0)}$ and $t\mapsto k_t\in\RR^{(d+d_0)}$. We  compute $d\YY_t$ applying It\^o's formula to this ansatz, and using the expression for $d\XX_t$ given by the forward equation. Identifying term by term the result with the right hand side of the backward component of the above FBSDE  we obtain the following system of ordinary differential equations:
\begin{equation}
\label{fo:first_system}
\begin{cases}
&\hskip -10pt
0=\dot{K}_t  - \frac12 K_t \BB_0 R_0^{-1} \BB_0^\dagger K_t + K_t \LL_0^{(cl)}(t) + \LL_0^{(cl)}(t)^\dagger K_t+ \FF_0,\quad K_T = 0\\
&\hskip -10pt
0=\dot{k}_t  +\bigl( \LL_0^{(cl)}(t)^\dagger -\frac12  K_t \BB_0 R_0^{-1} \BB_0^\dagger\bigr) k_t + K_t\CC_0^{(cl)}(t) + 2 f_0 ,\;\;\;k_T =0.
\end{cases}
\end{equation}
For any choice of a continuous strategy $t\mapsto (\phi_0(t),\phi_1(t),\phi_2(t),\phi_3(t))$, the first equation is a standard matrix Riccati differential equation. Since the coefficients are continuous and $\FF_0$ is positive definite, the equation admits a unique global solution over $[0,T]$ for any $T>0$. Recall that $R_0$ is symmetric and positive definite. Injecting the solution $t\mapsto K_t$ into the second equation yields a linear ordinary differential equation with continuous coefficients for which the global unique solvability also holds. Therefore the FBSDE \eqref{fo:fbsde_xx_yy} is uniquely solvable and the optimal control exists and is given by:
\begin{equation}
\label{fo:mmmfg_cl_major_opt}
\alpha_t^{0 *} = - \frac12 R_0^{-1}\BB_0^\dagger K_t \XX_t - \frac12 R_0^{-1}\BB_0^\dagger k_t,
\end{equation}
which is an affine function of $X_t^0$ and $\bar X_t$. 

\vskip 6pt
We now turn to representative minor player optimization problem. We assume that the major player uses the feedback strategy $\alpha^0_t=\phi^0(t, X_t^0, \bar X_t)$ and the representative of the other minor players uses the feedback strategy $\alpha_t=\phi(t, X_t, X_t^0, \bar X_t)$ of the forms \eqref{fo:alpha0} and \eqref{fo:alpha} respectively. These choices lead to the dynamics of the state $\XX_t = [\overline X_t^\dagger, X_t^{0\dagger}]^\dagger$ given by:
$$
d\XX_t=[\LL^{(cl)}(t) \XX_t + \CC^{(cl)}(t) ]dt + \mathbb{D}_0 dW_t^0
$$
with:
\begin{equation*}
\LL^{(cl)}(t) =\left[\begin{array}{cc}L+F+B(\phi_1(t) + \phi_3(t)) & G + B\phi_2(t)  \\ F_0 +  B_0\phi^0_2(t)  & L_0  + B_0\phi^0_1(t) \end{array}\right], 
\quad
\CC^{(cl)}(t) =\left[\begin{array}{c} B\phi_0(t) \\ B_0\phi_0^0(t) \end{array}\right].
\end{equation*}
We wrote $\LL^{(cl)}(t)$ and $\CC^{(cl)}(t)$ instead of $\LL^{(cl),\phi^0,\phi}(t)$ and $\CC^{(cl),\phi^0,\phi}(t)$ in order to simplify the notation. In this environment, we search for the best response of a representative minor player trying to minimize as earlier,
$$
\inf_{\tilde\balpha\in\AA}\EE \left[\int_0^T [(\tilde X_t - [H_1, H] \XX_t - \eta)^{\dagger} Q(\tilde X_t - [H_1, H] \XX_t - \eta)+\tilde\alpha_t^\dagger R \tilde\alpha_t ] dt \right],
$$
where the dynamics of the controlled state $\tilde X_t$ are given as before by:
$$
d\tilde X_t = (L \tilde X_t + B \tilde\alpha_t + [F,G]\XX_t) dt + D dW_t.
$$
Again the process $\mathbb{X}_t$ is merely part of the random coefficients of the optimization problem. We introduce the reduced Hamiltonian:
\begin{equation*}
\begin{split}
&H^{(r),\phi^0,\phi}(t,\tilde x,\tilde y,\tilde\alpha) = 
\tilde y^{\dagger}(L \tilde x + B \tilde\alpha + [F,G] \mathbb{X}_t)\\ 
&\hskip 40pt
+ (\tilde x - [H_1, H] \mathbb{X}_t - \eta)^{\dagger} Q(\tilde x - [H_1,H] \mathbb{X}_t - \eta)+\tilde\alpha^\dagger R \tilde\alpha.
\end{split}
\end{equation*}
and we find that the optimal control is given by $\tilde\alpha_t^*=-\frac12 R^{-1}B^\dagger Y_t$, where $(\tilde X_t, \XX_t,\tilde Y_t, \tilde Z_t,\tilde Z^0_t)_{0\le t\le T}$ solves the linear FBSDE:
\begin{equation*}
\begin{cases}
&\hskip -10pt
d\tilde X_t = (L\tilde X_t -B R^{-1}B^{\dagger}\tilde Y_t +  [F,G]\XX_t) dt + D dW_t\\
&\hskip -10pt
d\XX_t=[\LL^{(cl)}(t) \XX_t + \CC^{(cl)}(t) ]dt + \mathbb{D}_0 dW_t^0\\
&\hskip -10pt
d\tilde Y_t = -(L^{\dagger} \tilde Y_t +Q \tilde X_t - Q [H_1,H] \mathbb{X}_t - Q\eta) dt + \tilde Z_t dW_t + \tilde Z_t^0 dW_t^0,\;\;Y_T =0.
\end{cases}
\end{equation*}
Again we search for a solution of the form $\tilde Y_t = \SS_t \XX_t + S_t \tilde X_t + s_t$ for continuous deterministic functions $t\mapsto \SS_t\in\RR^{d\times (d+d_0)}$, $t\mapsto S_t\in\RR^{d\times d}$ and $t\mapsto s_t\in\RR^d$. Proceeding as before, we see that these functions provide a solution to the above FBSDE if and only if they solve the system of ordinary differential equations:
\begin{equation}
\label{fo:second_system}
\begin{cases}
&\hskip -10pt
0 =\dot{S}_t + S_t L  + L^\dagger S_t - S_tBR^{-1}B^\dagger S_t + Q,\;\;\;\;S_T = 0\\
&\hskip -10pt
0 =\dot{\SS}_t+ \SS_t \LL^{(cl)}(t) + L^\dagger  \SS_t - S_tBR^{-1}B^\dagger  \SS_t + S_t [F,G] - Q[H_1,H],\;\;\; \SS_T = 0\\
&\hskip -10pt
0 = \dot{s}_t + (L^\dagger - S_tBR^{-1}B^\dagger) s_t + \SS_t \CC^{(cl)}(t) - Q\eta,\;\;\;\;s_T = 0.
\end{cases}
\end{equation}
The first equation is a standard symmetric matrix Riccati equation. As before, the fact that $Q$ is symmetric and non-negative definite and $R$ is symmetric and positive definite imply that this Riccati equation has a unique solution on $[0,T]$. Note that its solution $S_t$ is symmetric and independent of the inputs feedback functions $\phi^0$ and $\phi$ giving the controls chosen by the major player and the other minor players. Injecting the solution $S_t$ into the second and third equations, leads to a linear system of ordinary differential equations which can be readily solved.
Given such a solution we find that the optimal control can be expressed as:
\begin{equation}
\label{fo:mmmfg_cl_minor_opt}
\tilde\alpha_t^* = -\frac12 R^{-1}B^\dagger[\SS_t \XX_t + S_t X_t + s_t]
\end{equation}
which is indeed an affine function of $X_t$, $X_t^0$ and $\bar X_t$.

\vskip 6pt
Now that the two optimization problems are solved, we can tackle the fixed point step. We just proved that the best response map leaves the set of affine controls of the forms \eqref{fo:alpha0} and \eqref{fo:alpha}  invariant. This suggests that we can look for a fixed point in this set. For such a fixed point, we must have:
$$
\alpha_t^{0,*} = \phi^0(t, X_t^0, \overline X_t)=\phi_0^0(t) +\phi_1^0(t)X^0_t +\phi^0_2(t)\overline X_t,
$$
and:
$$
\tilde\alpha_t^* = \phi(t, X_t, X_t^0, \bar X_t)=\phi_0(t) +\phi_1(t)X_t +\phi_2(t)X_t^0 +\phi_3(t)\overline X_t,
$$
which translates into the following equations:
\begin{align*}
[\phi^0_2(t),\phi^0_1(t)] = -\frac12 R_0^{-1}\mathbb{B}_0^\dagger K_t, &\;\;\;\;\; \phi^0_0(t) = -\frac12 R_0^{-1}\mathbb{B}_0^\dagger k_t,\\
[\phi_3(t),\phi_2(t)] = -\frac12 R^{-1}B^\dagger\SS_t, &\;\;\;\;\; \phi_1(t) = -\frac12 R^{-1}B^\dagger S_t,&\phi_0(t) = -\frac12 R^{-1}B^\dagger  s_t.
\end{align*}
To complete the construction of the equilibrium, it thus remain to determine the quantities $K_t$, $k_t$, $\SS_t$, $S_t$ and $s_t$ from the systems
\eqref{fo:first_system} and \eqref{fo:second_system}.
As we already noticed, the second equation of  \eqref{fo:first_system} can be used to determine $k_t$ from $K_t$. As for \eqref{fo:second_system}, $S_t$ can be obtained by solving the first equation on its own, and once this is done the third equation of  \eqref{fo:second_system} can be used to determine $s_t$ from $\SS_t$. In other words,  we can solve for $S_t$ by solving the first equation of \eqref{fo:second_system}, and then group the remaining four equations into two systems of ordinary differential equations as follows:
\begin{equation}
\label{fo:Riccatis}
\begin{cases}
&0 =\dot{K}_t + K_t[\LL(t) - \BB R^{-1}B^\dagger \SS_t] + [\LL(t) - \BB R^{-1}B^\dagger \SS_t]^\dagger K_t\\ 
&\hskip 165pt
-K_t \BB_0R_0^{-1}\BB_0^\dagger K_t + \LL_0\\
&0 =\dot{\SS}(t) + \SS_t\AA(t) + [L^\dagger - S_t BR^{-1}B^\dagger] \SS_t-  \SS_t\BB R^{-1}B^\dagger  \SS_t\\
&\hskip 75pt
 -  \SS_t \BB_0R_0^{-1}\BB_0^\dagger K_t + [S_t F - QH_1, S_tG - QH]\\
\end{cases}
\end{equation}
and
\begin{equation}
\label{fo:non_Riccatis}
\begin{cases}
&0=\dot k_t + [\LL(t) - \BB R^{-1}B^\dagger \SS_t]^\dagger k_t - K_t\BB_0R_0^{-1}\BB_0^\dagger k_t
- K_t\BB R^{-1}B^\dagger s_t + f_0\\
&0=\dot s_t + [L^\dagger - S_t BR^{-1}B^\dagger] s_t - \SS_t\BB_0R_0^{-1}\BB_0^\dagger k_t 
- \SS_t\BB R^{-1}B^\dagger s_t - Q\eta
\end{cases}
\end{equation}
with $0$ as terminal condition, where we used the notation:
$$
\LL(t) := \LL_0 -\left[ \begin{array}{cc}BR^{-1}B^\dagger S_t&0\\0&0\end{array}\right].
$$
The first system \eqref{fo:Riccatis} comprises two mildly coupled matrix Riccati equations, while the system \eqref{fo:non_Riccatis}, once the solutions of the
first system are identified and substituted for, is a plain linear system whose solution is standard.
In other words, the functions $t\mapsto k_t$ and $t\mapsto s_t$ can easily be determined once a solution $t\mapsto (K_t,\SS_t)$ of system \eqref{fo:Riccatis} is found. In essence, we proved the following verification theorem.

\begin{theorem}
If the system \eqref{fo:Riccatis} of matrix Riccati equations is well posed, then there exists a Nash equilibrium in the family of linear closed loop feedback controls, the optimal controls for the major and minor players being given by the strategies \eqref{fo:mmmfg_cl_major_opt} and \eqref{fo:mmmfg_cl_minor_opt}.
\end{theorem}

\section{Application}
\label{se:application}
In this final section, we apply the theoretical results derived above to a model of flocking inspired by the mean field game formulation proposed in \cite{NourianCainesMalhame} to generalize a basic descriptive model originally proposed by Cucker and Smale in \cite{CuckerSmale}.
In this section, we borrow from the terminology used in the dynamical systems literature on large population behavior, and we call the major player the \emph{leader} while the minor players are call \emph{followers}. However, the reader should not be misled by this terminology: we are not solving a leader-follower game, we are solving for a Nash equilibrium for the mean field game with major and minor players.

, in which the dynamics of a large population of agents are governed by forces depicting the mean reversion of individual's velocity to the mean velocity of the population. Later on, Huang (reference) formulates the flocking model into a mean field game, where the emergent behavior is obtained by the Nash equilibrium of the game. While early models of flocking does not involve any form of central coordination, several authors recently propose generalization of the flocking model by introducing leaders in the population. Such leader has a pivotal impact on the rest of the population. In this spirit, we generalize Huang's formulation of flocking mean field game by introducing a free-will leader pursuing a prescribed schedule of velocity. 

\vskip 4pt
Given a population of $N$ minor players (followers), we denote by $V_t^{0,N}$ the velocity of the major player (leader) at time $t$, and by $V_t^{n,N}$ the velocity of the $n$-th follower. The leader and the followers control the drifts of their velocities whose dynamics are
given as It\^o processes:

\begin{equation}
\label{fo:flocking_dynamics}
\begin{cases}
&dV_t^{0,N}=\alpha^0_t dt + \Sigma_0 dW_t^0\\
&dV_t^{n,N}=\alpha^n_t dt + \Sigma dW_t^n
\end{cases}
\end{equation}
where the $d$-dimensional Wiener processes $\{\bW^i=(W^i_t)_{0\le t\le T};\;i=0,1,\cdots,N\}$ are independent, and $\Sigma_0$ and $\Sigma$ are constant $d\times d$ matrices. We also assume that we are given a deterministic function $[0,T]\ni t \rightarrow \nu_t\in\RR^d$ representing the leader's free will, namely the velocity the major player would like to have while keeping a reasonable distance from the pack. If we denote by $\bar V_t^{N}:= \frac{1}{N}\sum_{n=1}^N V_t^{n,N}$ the average velocity of the followers, the objective of the leader is to minimize its expected costs over the horizon $T$:
\[
J^0 = \mathbb{E}\Bigl[\int_0^T \bigl(\lambda_0\|V^{0,N}_t - \nu_t\|^2 + \lambda_1\|V^{0,N}_t - \bar V^{N}_t\|^2 + (1 - \lambda_0 - \lambda_1) \|\alpha^0_t\|^2 \bigr)dt\Bigr]
\]
where $\lambda_0$ and $\lambda_1$ are positive real numbers satisfying $\lambda_0 + \lambda_1 \le 1$. Similarly, each follower faces a tradeoff between keeping up with the leader and staying close to its peers. So the objective of the $n$-th  follower is to minimize:
\[
J^n = \mathbb{E}\Bigl[\int_0^T \bigl( l_0\|V^{n,N}_t - V^{0,N}_t\|^2 + l_1\|V^{n,N}_t - \bar V^{N}_t\|^2 + (1 - l_0 - l_1) \|\alpha^n_t\|^2 \bigr)dt\Bigr]
\]
where $l_0$ and $l_1$ are positive reals satisfying $l_0 + l_1 \le 1$. While the above model is clearly linear quadratic, it does not fit in the framework used in this paper. However, it is plain to remedy this problem by simply doubling the state variable. More specifically, we define $X_t^0 := [V_t^0, V_t^0]$, $X_t := [V_t, V_t]$ and $\bar X_t := [\bar V_t, \bar V_t]$ and we pose:
\begin{equation*}
\begin{split}
L_0 = L = F_0 = F = G = \left[\begin{array}{cc}0 & 0 \\ 0&0 \end{array}\right],\;\;\; B_0 = B =  \left[\begin{array}{c}I \\ I \end{array}\right],\;\;D_0= \left[\begin{array}{c}\Sigma_0 \\ \Sigma_0 \end{array}\right],\;\;D= \left[\begin{array}{c}\Sigma \\ \Sigma \end{array}\right]\\
H = \left[\begin{array}{cc}I & 0 \\ 0&0 \end{array}\right],\;\;H_0 = H_1 = \left[\begin{array}{cc}0 & 0 \\ 0&I \end{array}\right],\;\;Q_0 = \left[\begin{array}{cc}\lambda_0I & 0 \\ 0&\lambda_1I \end{array}\right],\;\;Q = \left[\begin{array}{cc}l_0I & 0 \\ 0&l_1I \end{array}\right]\\
\eta_0(t) = \left[\begin{array}{c}\nu(t) \\ 0 \end{array}\right],\;\;\eta = \left[\begin{array}{c}0 \\ 0 \end{array}\right],\;\;R_0 = (1 - \lambda_0 - \lambda_1)I,\;\;R = (1 - l_0 - l_1)I
\end{split}
\end{equation*}

We implemented the solution of this model in the $d=2$ dimensional case choosing 
$$
\nu(t) := [-2\pi \sin(2\pi t), 2\pi \cos(2\pi t)]
$$ 
for the leader's free-will. We also choose $\Sigma_0 = \Sigma = 0.5 I_2$. For a given choice of penalty coefficients $\lambda_0, \lambda_1, l_0, l_1$, we use Euler's method to solve numerically the system of matrix Riccati equation \eqref{fo:Riccatis} over the horizon $T = 5$, and computing closed loop Nash equilibrium strategies of for the leader and the representative follower in the mean field game limit.

We simulate the dynamics of the leader and $N$ followers defined in \eqref{fo:flocking_dynamics}, where we assign the equilibrium control strategies of the mean field game to the leader and each follower.

\begin{figure}[H]
\centering
\includegraphics[scale=0.42, trim = 3mm 0mm 5mm 0mm, clip=true]{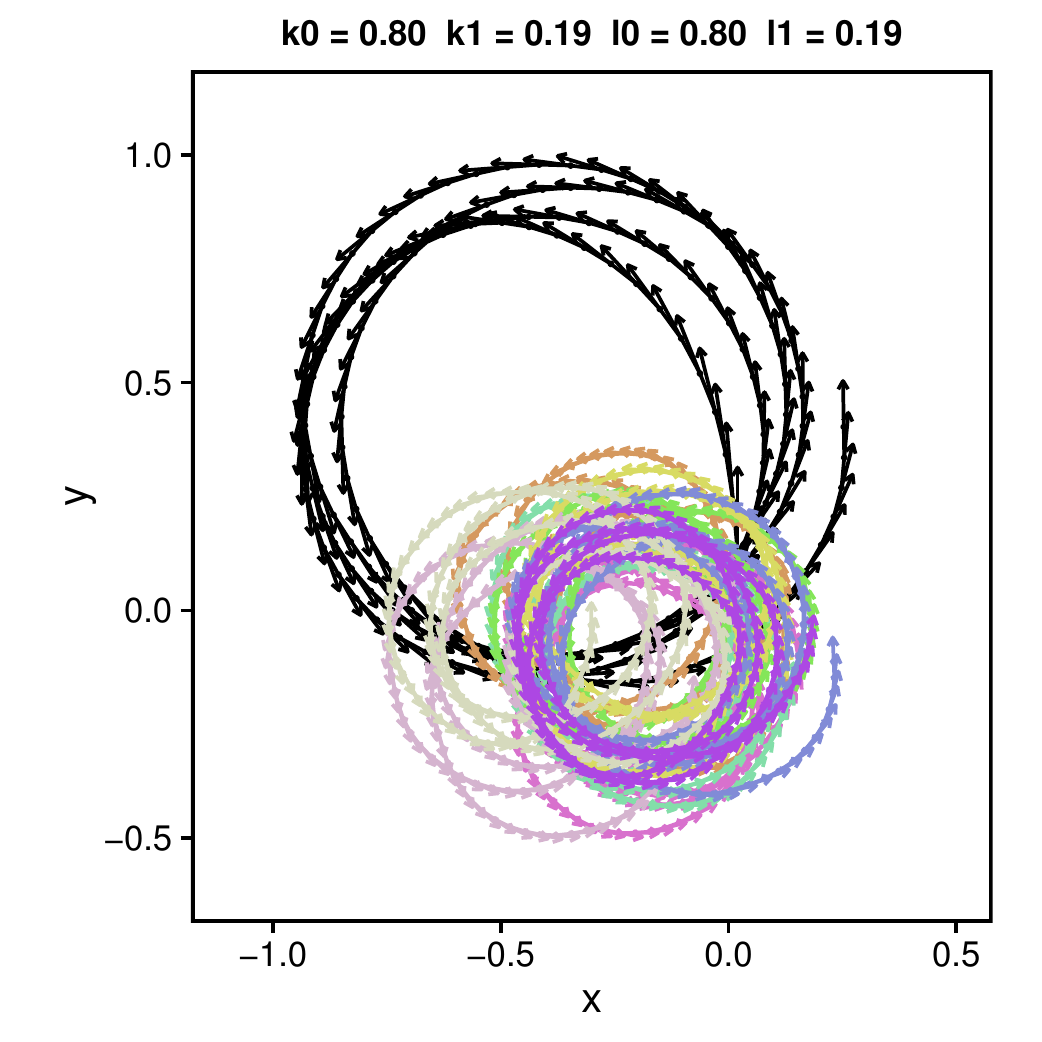}
\includegraphics[scale=0.42, trim = 3mm 0mm 3mm 0mm, clip=true]{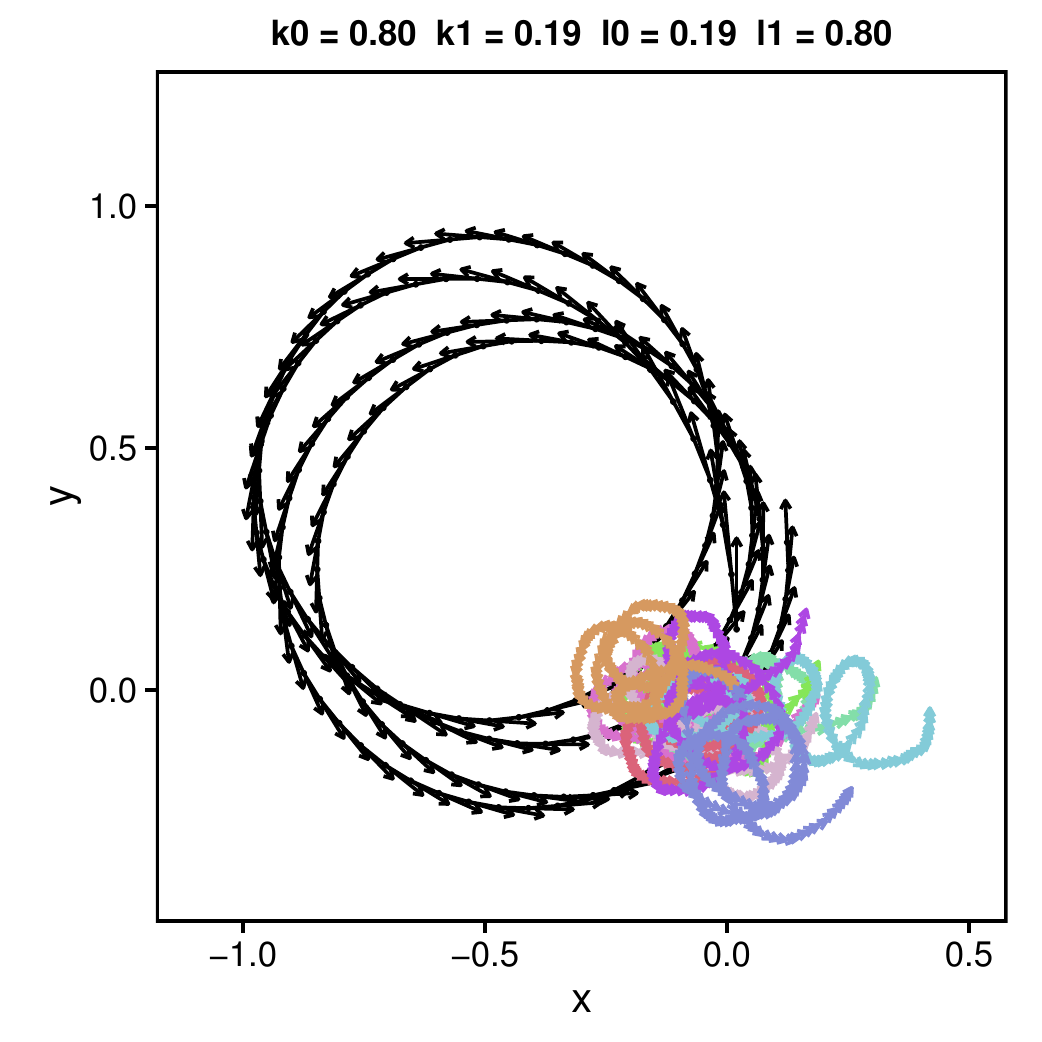}
\includegraphics[scale=0.42, trim = 3mm 0mm 3mm 0mm, clip=true]{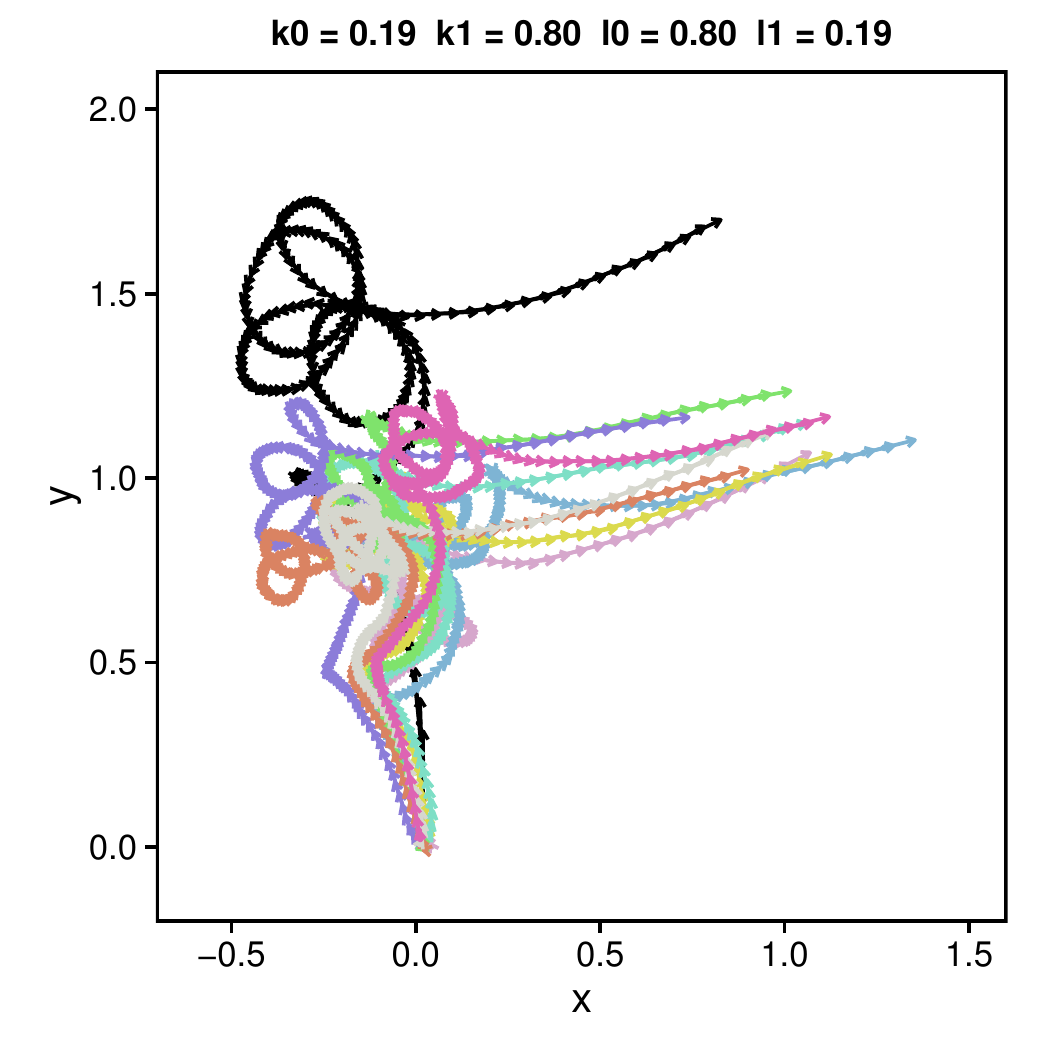}
\includegraphics[scale=0.42, trim = 3mm 0mm 3mm 0mm, clip=true]{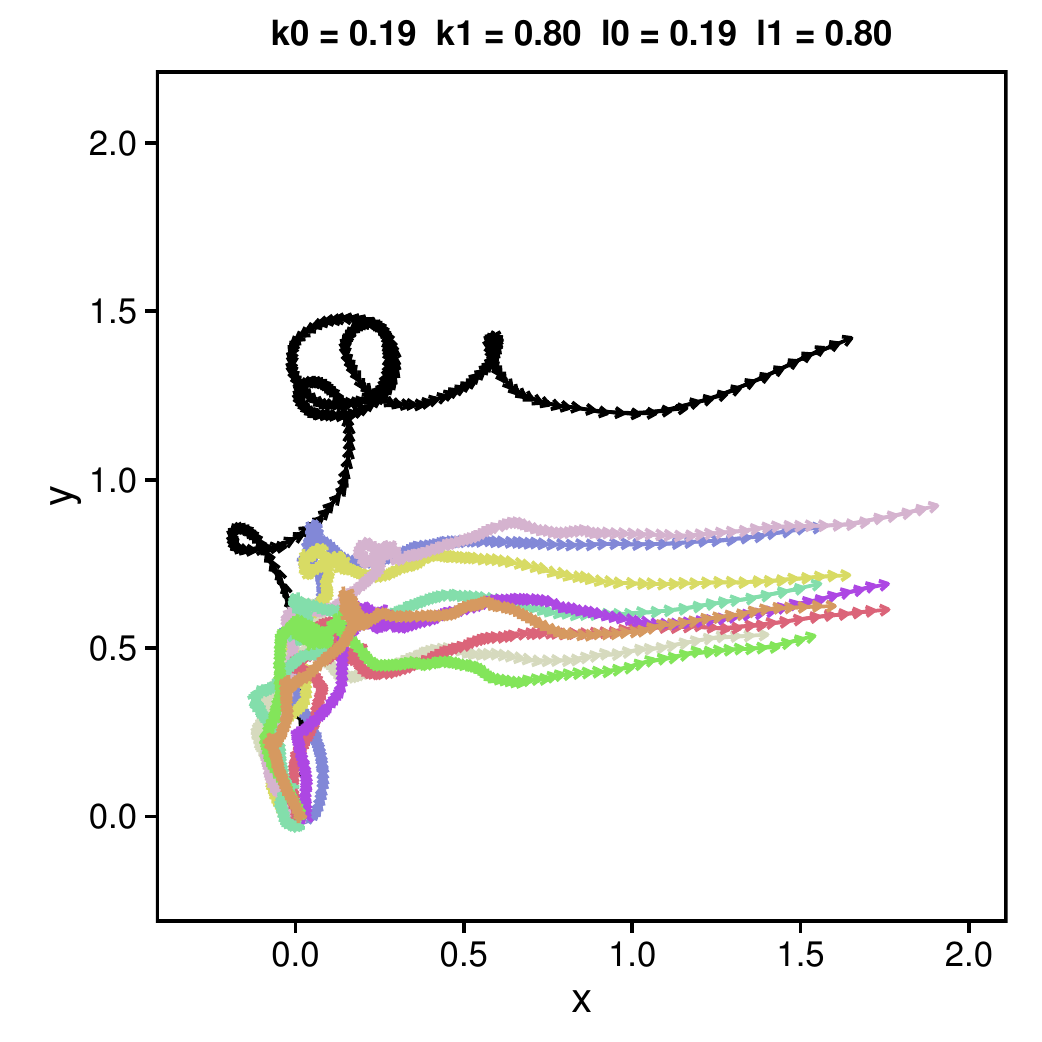}
\caption{Optimal velocity and trajectory of follower and leaders}
\label{figure:optimal_trajectory}
\end{figure}

\begin{figure}[H]
\centering
\includegraphics[scale=0.3, trim = 15mm 0mm 15mm 0mm, clip=true]{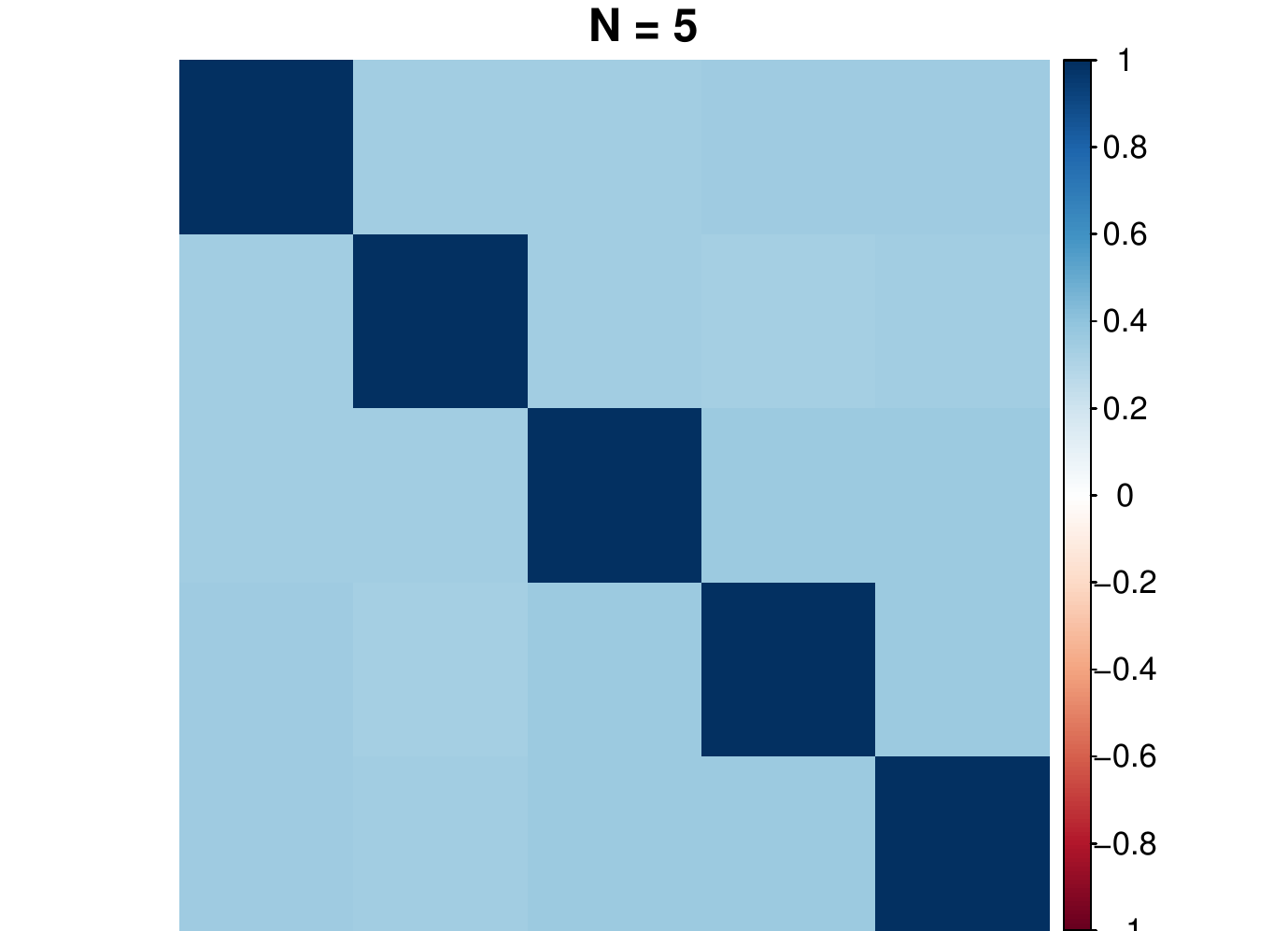}
\includegraphics[scale=0.3, trim = 15mm 0mm 15mm 0mm, clip=true]{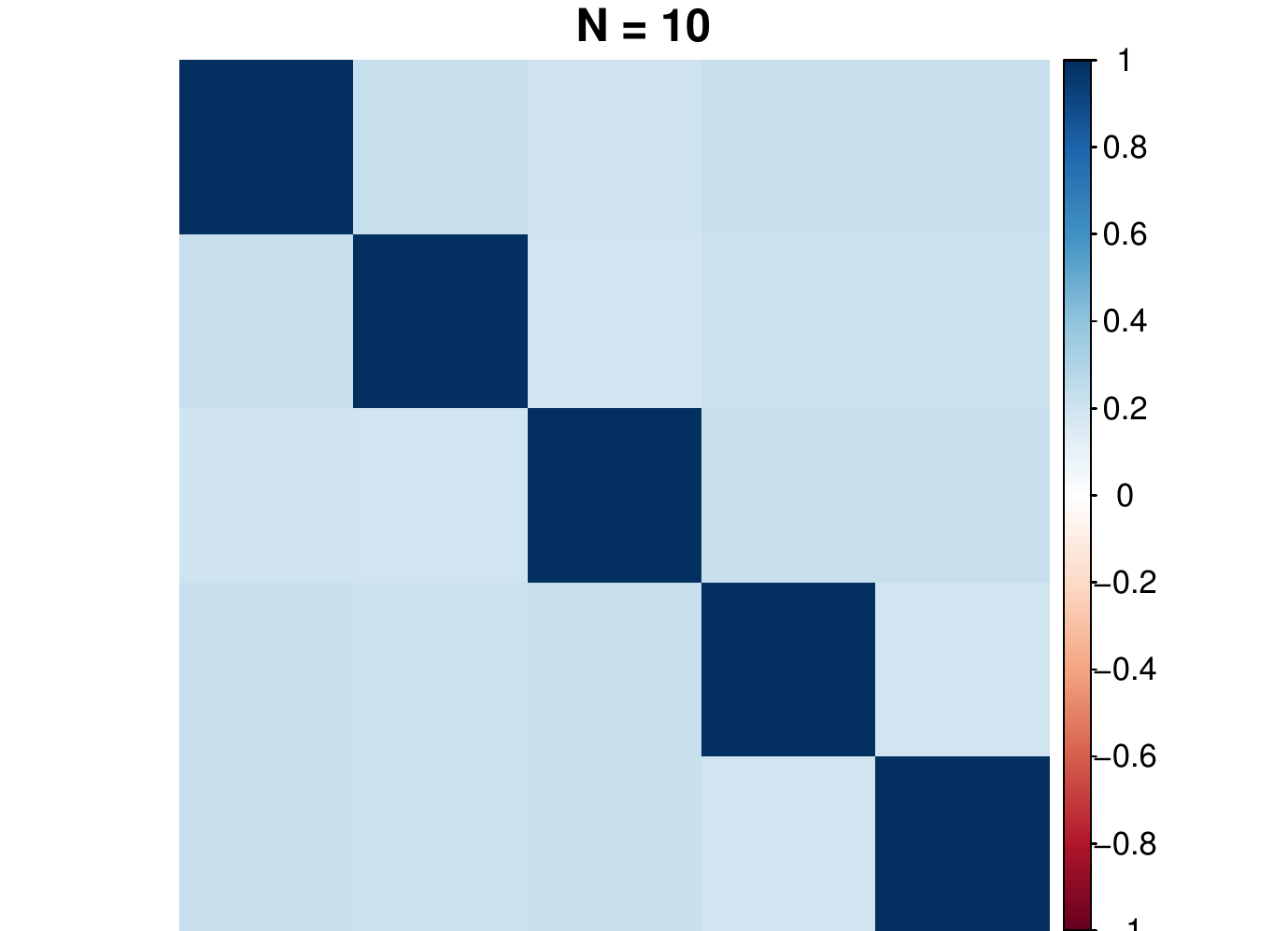}
\includegraphics[scale=0.3, trim = 15mm 0mm 15mm 0mm, clip=true]{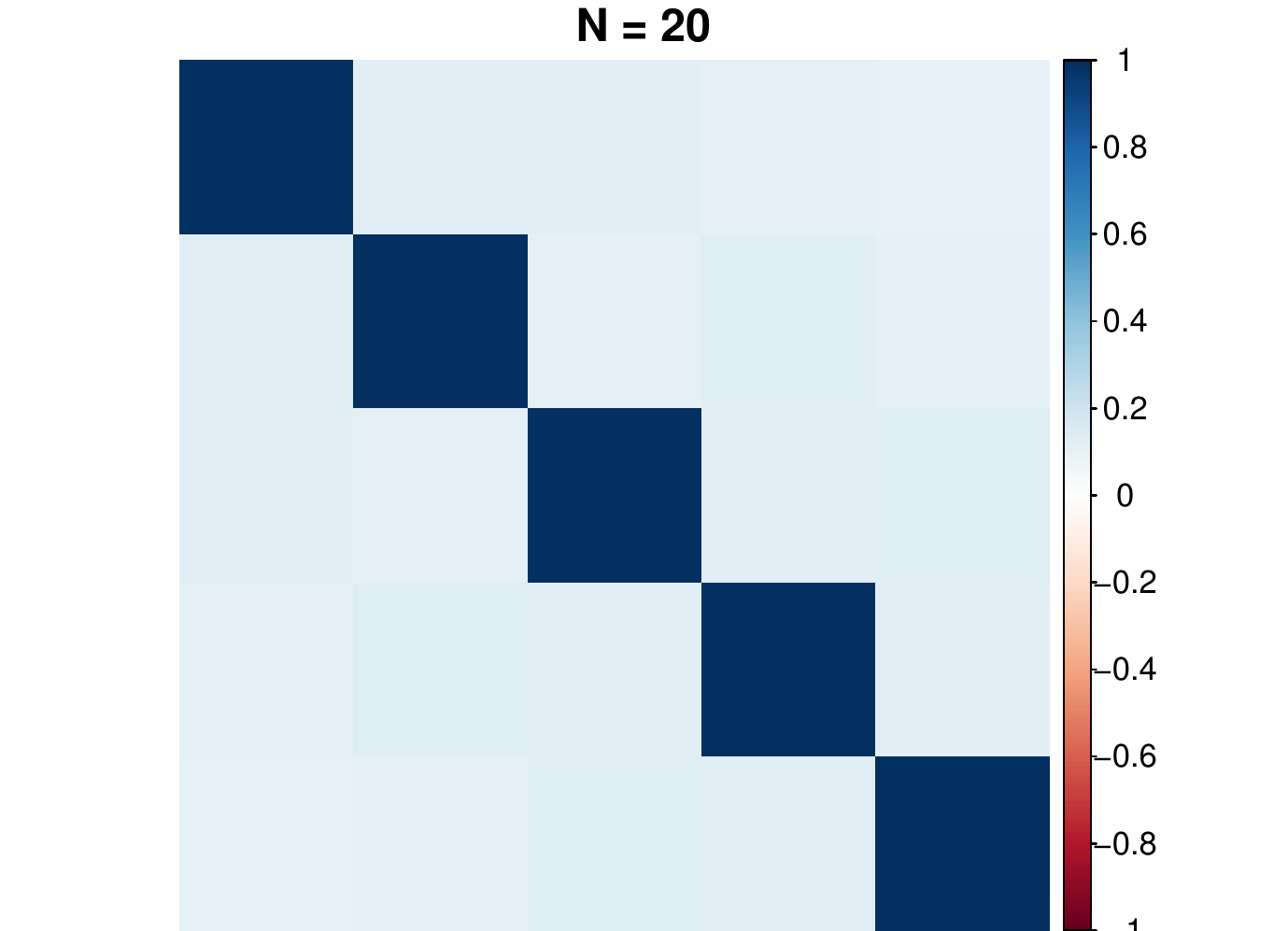}
\includegraphics[scale=0.3, trim = 15mm 0mm 15mm 0mm, clip=true]{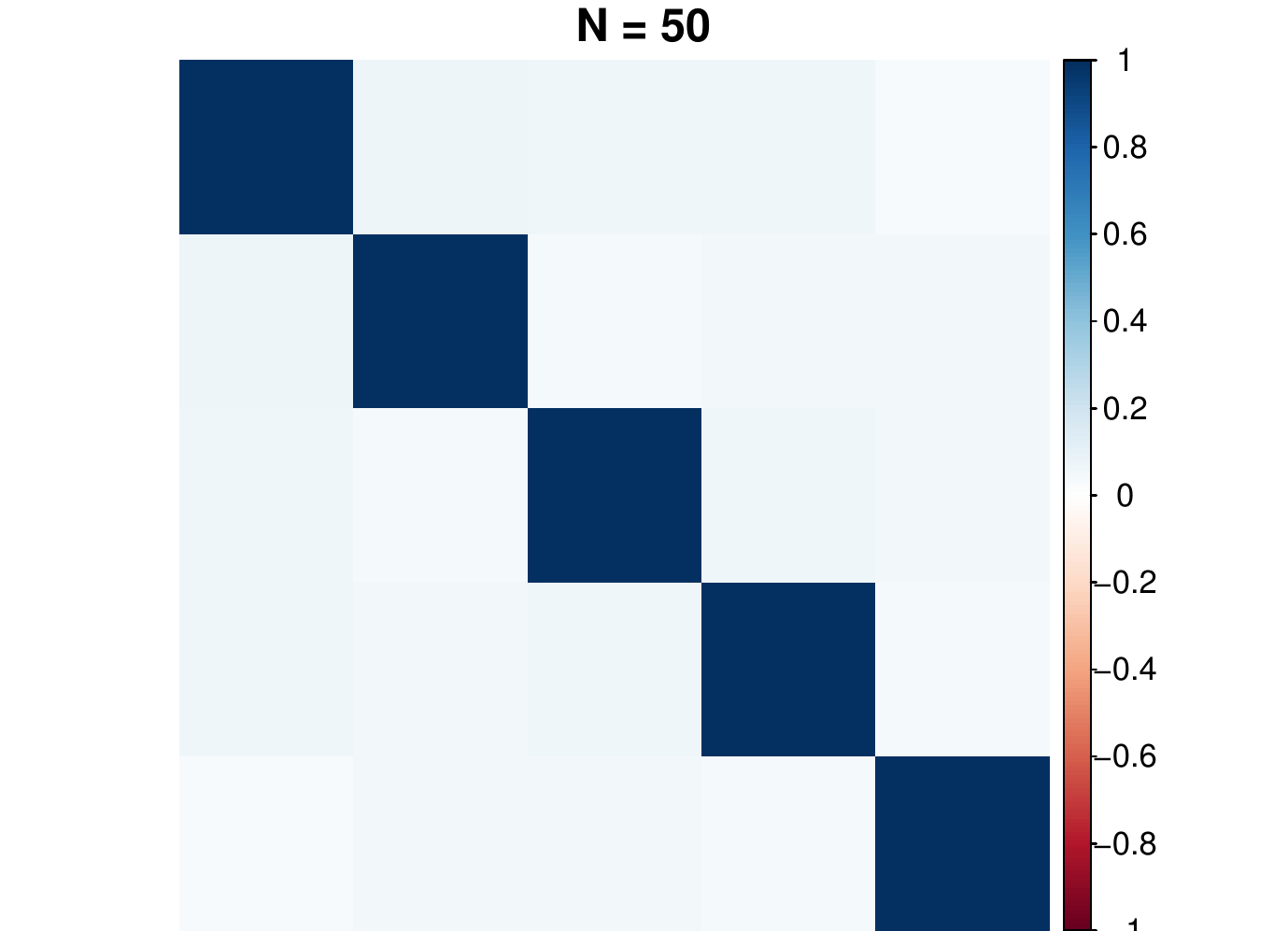}
\includegraphics[scale=0.3, trim = 15mm 0mm 3mm 0mm, clip=true]{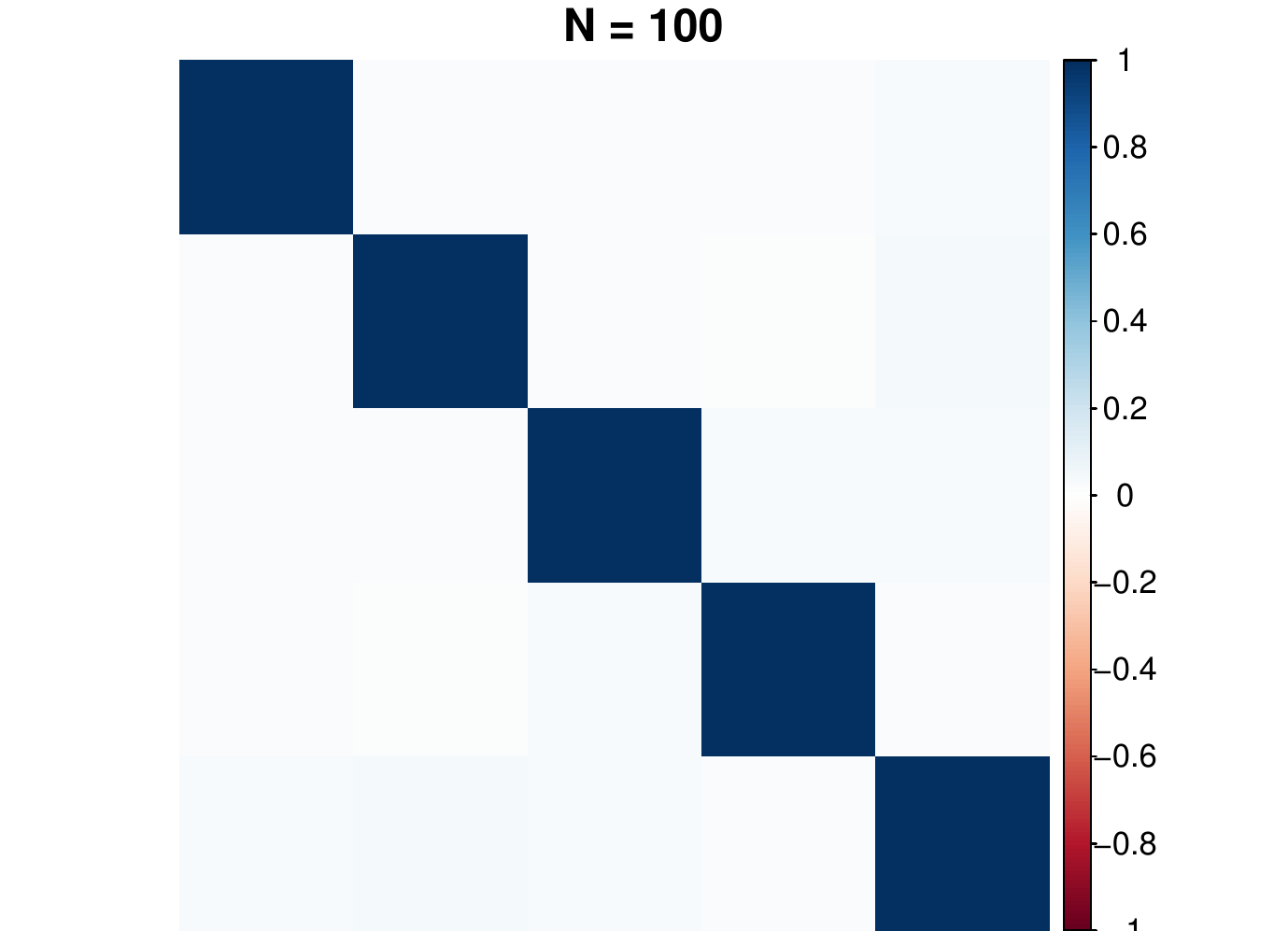}
\caption{Conditional correlation of followers' velocities}
\label{figure:propagation_chaos}
\end{figure}

\vskip 4pt
Figure \ref{figure:optimal_trajectory} shows the trajectories (points in the plane) and the velocities (arrows) of the flock. The leader's trajectory is plotted in black and those of the followers in color. We observe that the prescribed velocity $\nu$ is best followed by the flock when the leader cares more about pursuing its objective and the followers are more committed to follow the leader, rather than sticking with the average of the population. Conversely, if the individuals attribute more importance to staying close with the population, the flock follows an erratic trajectory in the beginning and eventually reaches a common direction of movement.
Our simulation also gives a peak into the effect of propagation of chaos, which states that in the limit of an infinite number of followers, the velocities of the followers become independent conditioned on the shock process driving the leader's velocity. To visualize such an effect, for a given number of followers, say $N$, we fix a realization of the Wiener processes $\bW^0$ driving the dynamics of the leader's velocity. We simulate $S$ copies of the optimal paths $V_t^{0,N}$ and $V_t^{n,N}, n = 1,\dots,N$ where for each sample path we use the same Wiener process we fixed before for the leader, but independent copy of Wiener process for each of the followers. Then  for a given $t$, we compute the sample correlation matrix of $V_t^{i,N,(1)}, i = 1,\dots, 5$, which are the first components of the velocity of the first 5 followers at time $t$. Finally, we compute the average of the correlation matrix across time $t\le T$. Figure \ref{figure:propagation_chaos} displays the average correlation matrices for flocks of sizes $N=5, 10, 20, 50, 100$ obtained by following the procedure described above. It can be seen that the correlation between the followers' velocities dramatically reduces to 0 as the size of the flock grows. Indeed, the linearity of the leader and follower strategies implies that the whole system evolves as a vector-valued OU process, and the velocity of any individual at a given time is Gaussian. Since independence is equivalent to null correlation for Gaussian vectors, the convergence of the correlation matrices provides a strong evidence of the conditional propagation of chaos.

\bibliography{games}

\begin{thebibliography}{10}

\bibitem{BensoussanChauYam}
{\sc A.~Bensoussan, M.~Chau, and S.~Yam}, {\em Mean field games with a
  dominating player}, tech. rep., 2013.

\bibitem{CarmonaDelarue_book_I}
{\sc R.~Carmona and F.~Delarue}, {\em Probabilistic Theory of Mean Field Games:
  vol. I, Mean Field FBSDEs, Control, and Games}, Stochastic Analysis and
  Applications, Springer Verlag, 2017.

\bibitem{CarmonaDelarue_book_II}
\leavevmode\vrule height 2pt depth -1.6pt width 23pt, {\em Probabilistic Theory
  of Mean Field Games: vol. II, Mean Field Games with Common Noise and Master
  Equations}, Stochastic Analysis and Applications, Springer Verlag, 2017.

\bibitem{CarmonaZhu}
{\sc R.~Carmona and G.~Zhu}, {\em A probabilistic approach to mean field games
  with major and minor players}, Annals of Applied Probability, 26 (2014),
  pp.~1535--1580.

\bibitem{CuckerSmale}
{\sc F.~Cucker and S.~Smale}, {\em Emergent behavior in flocks}, {IEEE}
  Transactions on Automatic Control, 52 (2007), pp.~852--862.

\bibitem{Huang}
{\sc M.~Huang}, {\em Large-population lqg games involving a major player: the
  nash equivalence principle}, {SIAM} Journal on Control and Optimization, 48
  (2010), pp.~3318--3353.

\bibitem{JaimungalNourian}
{\sc S.~Jaimungal and M.~Nourian}, {\em Mean-field game strategies for a
  major-minor agent optimal execution problem}, tech. rep., University of
  Toronto, March 15, 2015.

\bibitem{NourianCainesMalhame}
{\sc M.Nourian, P.~Caines, and R.~Malham{\'e}}, {\em Mean field analysis of
  controlled {C}ucker-{S}male type flocking: Linear analysis and perturbation
  equations}, in Proceedings of the 18th {IFAC} World Congress, {M}ilan,
  {A}ugust 2011, 2011, pp.~4471--4476.

\bibitem{NguyenHuang1}
{\sc S.~Nguyen and M.~Huang}, {\em Linear-quadratic-{G}aussian mixed games with
  continuum-parametrized minor players}, SIAM Journal on Control and
  Optimization,  (2012).

\bibitem{NguyenHuang2}
\leavevmode\vrule height 2pt depth -1.6pt width 23pt, {\em Mean field {LQG}
  games with mass behavior responsive to a major player}, in 51th {IEEE}
  Conference on Decision and Control, 2012.

\bibitem{NourianCaines}
{\sc M.~Nourian and P.~Caines}, {\em $\epsilon$-nash mean field game theory for
  nonlinear stochastic dynamical systems with major and minor agents}, tech.
  rep., 2013.

\end{thebibliography}

\end{document}